\documentclass[12pt]{article}
\usepackage{amsmath}
\usepackage{amssymb}
\usepackage{amsthm}
\usepackage[square,sort,comma,numbers]{natbib}
\usepackage{bbm, dsfont}
\usepackage{algorithmic}
\usepackage{algorithm}
\usepackage{caption}
\usepackage{subcaption}
\usepackage{graphicx}
\usepackage{epstopdf}
\usepackage{tabularx}
\usepackage{rotating}
\usepackage{multirow}
\usepackage{color, colortbl}

\newcommand{\Nens}{N} 
\newcommand{\Nobs}{m} 
\newcommand{\Nstate}{n} 
\newcommand{\X}{{\bf X}} 
\newcommand{\x}{{\bf x}} 
\newcommand{\J}{\mathcal{J}} 
\newcommand{\lp}{\left (} 
\newcommand{\rp}{\right )} 
\newcommand{\lb}{\left [} 
\newcommand{\rb}{\right ]} 
\renewcommand{\ln}{\left \|} 
\newcommand{\rn}{\right \|}
\renewcommand{\ll}{\left \{} 
\newcommand{\rl}{\right \}}
\newcommand{\B}{{\bf B}} 
\newcommand{\R}{{\bf R}} 
\newcommand{\N}{M} 
\newcommand{\y}{{\bf y}} 
\renewcommand{\H}{{\bf H}} 
\newcommand{\xm}{{\overline{\bf x}}} 
\newcommand{\xb}{{\bf x}^{\rm b}} 
\newcommand{\w}{\alpha} 
\newcommand{\W}{{\boldsymbol \alpha}} 
\newcommand{\I}{{\bf I}} 
\newcommand{\M}{\mathcal{M}} 
\newcommand{\Nor}{\mathcal{N}} 
\newcommand{\xt}{{\bf x}^{\rm true}} 
\newcommand{\Ho}{{\mathcal{H}}} 
\newcommand{\dx}{{{\boldsymbol \delta} {\bf x}}}
\newcommand{\DX}{{{\boldsymbol \delta} {\bf X}}}
\newcommand{\Q}{{\bf Q}} 
\renewcommand{\d}{{\bf d}} 
\newcommand{\radius}{\Delta} 
\newcommand{\m}{\mathcal{Q}} 
\newcommand{\sa}{{\bf s_{\W}}}
\newcommand{\JE}{\mathcal{J}_{\rm ens}} 
\newcommand{\La}{\mathcal{L}} 
\renewcommand{\P}{{\bf P}} 
\renewcommand{\Re}{\mathbbm{R}}
\renewcommand{\S}{{\bf S}} 
\newcommand{\errobs}{{\boldsymbol \epsilon}} 
\newcommand{\K}{{\bf K}} 

\newcommand{\DXS}{\delta{\bf \X}^{s}}
\newcommand{\nPOD}{{r}} 
\newcommand{\WPOD}{{\boldsymbol \beta}}
\newcommand{\wPOD}{{\beta}}

\newcommand{\g}{{\bf g}} 
\newcommand{\G}{{\bf G}} 
\newcommand{\dxt}{{\bf s}} 
\newcommand{\Z}{{\bf Z}} 
\newcommand{\PODB}{{\boldsymbol \Phi}} 
\newcolumntype{N}{>{\centering\arraybackslash} m{0.32\textwidth} }
\newcolumntype{V}{>{\centering\arraybackslash} m{0.005\textwidth} }

\newcommand{\bs}{{\rm b}} 
\newcommand{\as}{{\rm a}} 
\newcommand{\basis}{{\boldsymbol \Psi}}
\newcommand{\basisi}{{\boldsymbol \psi}}

\begin{document}

\thispagestyle{empty}
\setcounter{page}{0}

\begin{Huge}
\begin{center}
Computational Science Laboratory Technical Report CSL-TR-05-2014\\
\today
\end{center}
\end{Huge}
\vfil
\begin{huge}
\begin{center}
Elias D. Nino and Adrian Sandu
\end{center}
\end{huge}

\vfil
\begin{huge}
\begin{it}
\begin{center}
``A Derivative-Free Trust Region Framework for Variational Data Assimilation''
\end{center}
\end{it}
\end{huge}
\vfil

\begin{large}
\begin{center}
Computational Science Laboratory \\
Computer Science Department \\
Virginia Polytechnic Institute and State University \\
Blacksburg, VA 24060 \\
Phone: (540)-231-2193 \\
Fax: (540)-231-6075 \\ 
Email: {\rm enino@vt.edu},{\rm sandu@cs.vt.edu} \\
Web: {\rm http://csl.cs.vt.edu}
\end{center}
\end{large}

\vspace*{1cm}

\begin{center}
\begin{tabular}{c}
\includegraphics[width=1.5cm]{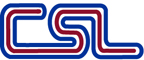}\\ 
.  \\
\includegraphics[width=1.5cm]{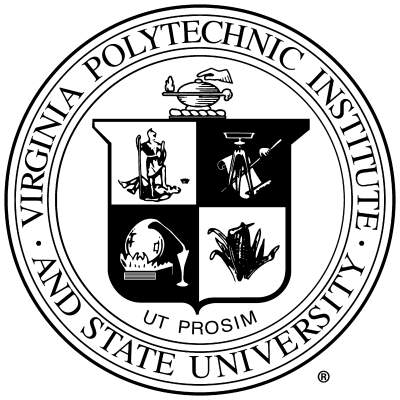} \\
\end{tabular}
\end{center}

\title{A Derivative-Free Trust Region Framework for Variational Data Assimilation}

\author{Elias D. Nino and Adrian Sandu \\
Computational Science Laboratory, Department of Computer Science \\
Virginia Polytechnic Institute and State University \\
Blacksburg, VA 24060, USA \\
enino@vt.edu, sandu@cs.vt.edu}

\maketitle
\tableofcontents

\begin{abstract}
This study develops a hybrid ensemble-variational approach for solving data assimilation problems. The method, called TR-4D-EnKF, is based on a trust region  framework and consists of three computational steps. First an ensemble of model runs is propagated forward in time and snapshots of the state are stored. Next, a sequence of basis vectors is built and a low-dimensional representation of the data assimilation system is obtained by projecting the model state  onto the space spanned by the ensemble perturbations. Finally, the low-dimensional optimization problem is solved in the reduced-space using a trust region approach; the size of the trust region is updated according to the relative decrease of the reduced order surrogate cost function. The analysis state is projected back onto the full space, and the process is repeated with the current analysis serving as a new background. A heuristic approach based on the trust region size is proposed in order to adjust the background error statistics from one iteration to the next. Experimental simulations are carried out {using the Atmospheric General Circulation Model (SPEEDY)}. The results show that TR-4D-EnKF is an efficient computational approach, and is more accurate than the current state of the art 4D-EnKF implementations such as the POD-4D-EnKF and the Iterative Subspace Minimization methods. 
\end{abstract}

\section{Introduction}
\label{sec:introduction}

Data assimilation \cite{Sandu_2011_assimilationOverview} is the process of estimating the true state $\xt_{\N} \in \Re^n$ of a dynamical system at the current time $t_{\N}$ given 
a history of prior evolution and noisy observations of the state at times $t_k$
\begin{eqnarray}
\label{eq:observations}
\y_k = \Ho_{k} \lp \xt_k \rp + \errobs_k \in \Re^{\Nobs \times 1}\,, \quad 0 \le k \le \N\,.
\end{eqnarray}
Here $\Nstate$ is the number of components in the model state, $\Nobs$ is the number of observed components from $\xt$, $\Ho_k:\Re^{\Nstate} \rightarrow \Re^\Nobs$ is the observation operator, $\errobs_k \in \Re^{\Nobs \times 1}$  is the error associated to the $k$-th observation time, and $\N$ is the number of observation times. Typically, observational errors are assumed to be normal distributed $ \errobs_k \sim \Nor({\bf 0}_{\Nobs}, \R_k)$ where ${\bf 0}_{\Nobs}$ is the $\Nobs$-th dimensional vector whose components are all zeros, and $\R_k \in \Re^{\Nobs \times \Nobs}$ is the data error covariance matrix at the assimilation time $t_k$. 

A dynamical model encapsulating our knowledge of the physical laws approximates the evolution of the dynamical system. 
The evolution of the model state $\x$ is given by 
\begin{eqnarray}
\label{eq:model-operator}
\x_{k+1} = \M_{t_k \rightarrow t_{k+1}} \lp \x_k \rp\,, \quad 0 \le k \le \N-1\,,
\end{eqnarray}
where $\M$ represents a nonlinear model solution operator (e.g., which simulates the evolution of the ocean or the atmosphere). 

Two families of methods, statistical filters and variational,  are widely used to solve data assimilation problems. Representative methods of those classes are the Ensemble Kalman Filter (EnKF) and the Four-Dimensional Variational Method (4D-Var), respectively. In EnKF an ensemble of model runs is propagated in time; when data is available the filtering step generates an {\it analysis ensemble} whose empirical mean is an estimator for $\xt$. Strong constraint 4D-Var seeks an {\it analysis initial state} such that the corresponding forecast best fits the observations within the assimilation window. It is well-accepted that both methods face specific challenges in practical applications where $\Nstate \sim 10^9$. For instance, ensemble-based filters suffer from statistical sampling errors, while variational methods require adjoint models which are labor-intensive to develop and computationally expensive to run.  

Hybrid methods have been proposed in order to combine the strengths of EnKF and 4D-Var methods. A decomposition of the background errors in  components that are analyzed and components that are ignored has been used to estimate posterior covariances \cite{Sandu_2010_hybridCovariance}, and the theoretical similarities between the two approaches have been used to construct look-ahead assimilation techniques  \cite{Sandu_2011_subspaceDA}. Other hybrid approaches are based on model reduction and/or space reduction \cite{Chen2011,Lucia2004,Song2013}.
A discussion of model reduction techniques is given in \cite{Olivier2009}.
In this paper, we focus on the reduced-space approach where a subspace of the state space is identified, the variational problem is solved in this subspace, then the analysis is projected back onto the model space. The new solution can be treated as a new background and the process is repeated. The reduced space data assimilation approaches available in the literature update the solution in the model space  unconditionally. No available method provides a relation between the analysis at the current iteration and its associated error statistics (i.e., the initial background error covariance matrix is assumed to hold for all the ensembles at all iterations). This fact is important in the sampling process, as the uncertainty associated with the analysis state decreases as the iterations progress. 

In this work we formulate a hybrid data assimilation algorithm in the context of derivative-free optimization. 
A rigorous  Trust Region (TR) framework is proposed where the TR radius in the model space is linked with the spread of the ensemble members and with the quality of the solutions found in the reduced-space. The new method is named TR-4D-EnkF.

The remainder of the paper is organized as follows. Section \ref{sec:preliminaries} reviews the current state of the art ensemble-based approaches to data assimilation. Section \ref{sec:proposed-method} develops the novel derivative free TR-4D-EnkF method. Numerical results using the Lorenz-96 and the quasi-geostrophic models are reported in Section \ref{sec:experimental}, and conclusions are presented in Section \ref{sec:conclusions}.
%

\section{Four-dimensional ensemble-based approaches to data assimilation}
\label{sec:preliminaries}
%
EnKF \cite{Evensen2009} is one of the most widely used methods in data assimilation due to its simple formulation and ease of implementation. Normality assumptions are made on both the background and data errors \cite{Evensen2009A}. The method contains two steps, the forecast and the analysis.

The {\it prior} (background) distribution is approximated 
by an ensemble of $\Nens+1$ model state samples
\begin{eqnarray}
\label{eq:initial-ensemble}
\X_0 = \lb \x_{0}^{\bs(1)},\,\x_{0}^{\bs(2)},\,\ldots,\,\x_{0}^{\bs(\Nens+1)}\rb \in \Re^{\Nstate \times (\Nens+1)} \,,
 \end{eqnarray}
with the empirical moments 
\begin{subequations}
\label{eq:empirical-moments}
\begin{equation}
\label{eq:ensemble-mean}
\displaystyle \xm_0 = \frac{1}{\Nens+1} \cdot \sum_{i=1}^{\Nens+1} \x_{0}^{\bs(i)} \in \Re^{\Nstate \times 1}\,,
\end{equation}
\begin{equation}
\label{eq:sample-covariance-mean}
\displaystyle \S_0 = \frac{1}{\Nens} \cdot \DX_0  \cdot {\DX_0}^T \in \Re^{\Nstate \times \Nstate}\,,
\end{equation}
\end{subequations}
where $\x_0^{\bs(i)}$ is the i-th ensemble member and the columns of matrix 
\begin{eqnarray}
\label{eq:anomalies-enkf}
\DX_0 = \lb \dx_0^{(1)}, \dx_0^{(2)}, \ldots, \dx_{0}^{(\Nens+1)}\rb \in \Re^{\Nstate \times (\Nens+1)}\,
\end{eqnarray}
are given by $\dx_0^{(i)} = \x_0^{\bs(i)}-\xm_0$, for $1 \le i \le \Nens+1$. Prior any measurement, the background state $\x_0^{\bs} \approx \xm_0$ provides the best estimation to $\xt_0$. 

In the forecast step the background  ensemble \eqref{eq:initial-ensemble} is obtained by an ensemble of model runs that propagate each model
state to the current time $t_{k}$. 

In the analysis step a {\it posterior} (analysis) ensemble is constructed by making use of the observation $\y_k$ and by applying the Kalman filter to each background ensemble member:
\begin{eqnarray}
\displaystyle \x_k^{\as(i)} = \x_k^{\bs(i)}+\K \lb {\y_k^{s(i)}} +  \errobs_k^{s(i)} - \H_k \cdot \x_k^{\bs(i)} \rb \,, \quad 1 \le i \le \Nens+1\,,
\end{eqnarray}
where $ \Ho_k' = \H_k \in \Re^{\Nobs \times \Nstate}$ is a linearized observation operator at time $t_k$, $\y_k^{s(i)} \sim \Nor(\y_k,\R_k)$ are the observations
$\y_k$ with added synthetic noise $\errobs_k^{s(i)}  \sim \Nor({\bf 0}_{\Nobs}\,,\R_k)$, and the Kalman gain matrix is $\K = \S_k \cdot \H_k^T \lb \H_k \cdot \S_k \cdot \H_k^T + \R_k \rb^{-1} \in \Re^{\Nstate \times \Nobs}$. 
The ensemble members are further propagated in time
\begin{eqnarray}
\label{eq:propagation}
\displaystyle \x^{\bs(i)}_{k+1} :=  \M_{t_k \rightarrow t_{k+1}} \lp \x^{\as(i)}_{k} \rp \,,
\end{eqnarray}
to obtain the background ensemble for the forecast step. 
EnKF can provide flow-dependent error estimates of the background errors (with the Monte Carlo methods) \cite{Nino2012,Nino2014}, but it does not have the ability to assimilate the observation data available at distributed times. 

4D-Var considers cost functions of the form
\begin{eqnarray} 
\label{eq:J4VAR}
\displaystyle \J(\x_0) &=& \underbrace{\frac{1}{2} \ln \x_0-\x_0^{\bs}\rn_{\B_0^{-1}}^2}_{\J^{\bs}(\x)}+ \underbrace{\frac{1}{2} \sum_{k=0}^{\N} \ln \y_{k} - \Ho \lp \x_k \rp \rn_{\R_k^{-1}}^2}_{\J^{\rm o}(\x)} \,,
\end{eqnarray}
where $\J^{\bs}(\x)$ and $\J^{\rm o}(\x)$ are known as the background and observation cost functions, respectively. 
The cost function \eqref{eq:J4VAR} is the negative logarithms of the a posteriori probability density when all the data and background
errors are normally distributed. 
The maximum likelihood estimate of the initial state is then obtained by minimizing the cost function,
i.e.,  the analysis step is computed by solving the optimization problem
\begin{eqnarray}
\label{eq:analysis-full-model}
\displaystyle \x_0^{\as} = \underset{\x_0} {\mathrm{arg\,min}} \,\J(\x_0) \qquad \textnormal{subject to  \eqref{eq:model-operator}}.
\end{eqnarray}
The formulation of \eqref{eq:J4VAR} allows 4D-Var to assimilate data which appears at different observation times. 

The computation of the gradient \eqref{eq:J4VAR} with respect to the control variable $\x_0 \in \Re^{\Nstate \times 1}$ requires one forward and one adjoint model integration.
The construction of an adjoint model for real, large forecast models is an extremely labor-intensive process.  In order to avoid the implementation
of adjoint models  four dimensional ensemble Kalman filter methods (4D-EnKF) \cite{Zhang2011} have been recently proposed. They naturally propagate
flow dependent background covariance matrices via ensembles \cite{Hamill2000,Lorenc2003,Song2013,Chung2013}. Numerical experiments show 
robust performance with a small number of ensemble members \cite{Nusrat2013,Hoteit2013}. Moreover, the solution \eqref{eq:analysis-full-model} can be treated as the new background state in \eqref{eq:J4VAR}, which provides a better solution \cite{Candiani2013}.

4D-EnKF based methods are defined as follows. The initial ensemble \eqref{eq:initial-ensemble} is propagated in time and $\N+1$ snapshots of each 
background ensemble member state at time moments $t_0,t_1,\dots,t_\N$ along the trajectory are stored
\begin{eqnarray}
\label{eq:snapshots}
\displaystyle \X^{s} = \begin{bmatrix}
\x_{0}^{b(1)} & \x_{0}^{b(2)} & \ldots & \x_{0}^{b(\Nens+1)} \\
\x_{1}^{b(1)} & \x_{1}^{b(2)} & \ldots & \x_{1}^{b(\Nens+1)} \\
\vdots & \vdots & \ddots & \vdots \\
\x_{\N}^{b(1)} & \x_{\N}^{b(2)} & \ldots & \x_{\N}^{b(\Nens+1)} \\
\end{bmatrix} \in \Re^{(\Nstate \cdot (\N+1)) \times (\Nens+1)} \,.
\end{eqnarray}
Each entry of the background ensemble matrix $\X^s$ is an $\Nstate$-dimensional vector $\x_{k}^{b(i)}$ which represents the state of ensemble member 
$i$ at time $t_k$. The $i$-th column of $\X^s$ contains all the snapshots of the $i$-th ensemble member, and the $k$-th row of blocks corresponds to all ensemble member states at $t_k$. 

Consider now a trajectory of the model. The state $\x_k$ at $t_k$ is approximated by a linear combination of the anomalies (deviations from the mean)
\begin{eqnarray}
\label{eq:solution-in-ensemble-space}
\displaystyle \x_k = \xm_k + \sum_{i=1}^{\Nens} \w_{i} \cdot \underbrace{\lp \x_k^{b(i)}-\xm_k \rp}_{\basisi_k^{(i)}} = \xm_k + \basis_k \cdot \W\,,
\end{eqnarray}
where 
\begin{eqnarray}
\label{eq:background-kth-iteration}
\displaystyle \xm_k &=& \frac{1}{\Nens+1} \cdot \sum_{i=1}^{\Nens+1} \x_{k}^{b(i)} \in \Re^{\Nstate \times 1}\,, \\
\label{eq:matrix-of-anomalies-k}
\displaystyle \basis_k &=&  \lb \basisi_k^{(1)}, \basisi_k^{(2)}, \ldots, \basisi_k^{(\Nens)} \rb \in \Re^{\Nstate \times \Nens}\,,
\end{eqnarray}
and the time-independent weight vector
\begin{eqnarray*}
\displaystyle \W = \lb \w_1, \w_2, \ldots, \w_{\Nens}\rb^T \in \Re^{\Nens \times 1} \,,
\end{eqnarray*}
contains the coordinates of $\x_k$ in the ensemble space.

By replacing \eqref{eq:solution-in-ensemble-space} in \eqref{eq:J4VAR} and linearizing the observation operator $\Ho_k \approx \H_k$, the 4D-Var cost function 
 \eqref{eq:J4VAR} can be written  in the ensemble space as follows:
\begin{eqnarray}
  \label{eq:J4VAR-alpha}
\displaystyle \JE(\W) &=& \frac{1}{2}\, \ln \d^{\rm b} - \basis_0 \cdot \W \rn_{\B_0^{-1}}^2+ \frac{1}{2} \, \sum_{k=0}^{\N} \ln \d_{k}^{\rm o}- \Q_k \cdot \W \rn_{\R_k^{-1}}^2
\end{eqnarray}
where $\d^{\rm b} = \x_0^{\bs}-\xm_0 \in \Re^{\Nstate \times 1}$ and $\d_{k}^{\rm o} = \y_{k}-\H_k \cdot \xm_k \in \Re^{\Nobs \times 1}$ are the innovation vectors on the background and observations, respectively, and $\Q_k = \H_k \cdot \basis_k \in \Re^{\Nobs \times \Nens}$. 

The optimal solution in the ensemble space
\begin{eqnarray}
\label{eq:solution-ensemble-space}
 \W^{*} = \underset{\W} {\mathrm{arg \,min}} \,\JE(\W) \in \Re^{\Nens \times 1}\,,
\end{eqnarray}
provides an approximation of the analysis trajectory started from \eqref{eq:analysis-full-model}  through the relation
\begin{eqnarray}
\label{eq:analysis-state-sol-ensemble}
\displaystyle \x_k^{\as} = \x_k^{\bs}+\basis_k \cdot \W^{*} \in \Re^{\Nstate \times 1}\,.
\end{eqnarray}
The derivatives of \eqref{eq:J4VAR-alpha} are 
\begin{subequations}
\label{eq:J4VAR-derivatives-ensemble}
\begin{eqnarray}
\label{eq:J4VAR-derivat}
\displaystyle \nabla_\W \JE(\W)&=& \lb \basis_0^T \cdot \B_0^{-1} \cdot \basis_0  +  \sum_{k=0}^{\N}  \Q_k^T \cdot \R_k \cdot \Q_k  \rb \cdot  \W \\
\nonumber
&-& \lb \basis_0^T \cdot \B_0^{-1} \cdot \d^{\bs} + \sum_{k=0}^{\N}  {\Q_k^T \cdot \R_k^{-1} \cdot \d_{k}} \rb   \in \Re^{\Nens \times 1} \,, \\
\label{eq:J4VAR-derivat2}
\displaystyle \nabla^2_{\W,\W} \JE(\W)&=& \basis_0^T \cdot \B_0^{-1} \cdot \basis_0  +  \sum_{k=0}^{\N}  \Q_k^T \cdot \R_k \cdot \Q_k \in \Re^{\Nens \times \Nens}\,,
\end{eqnarray}
\end{subequations}
and the solution of the quadratic minimization problem \eqref{eq:solution-ensemble-space} is
\begin{eqnarray}
\label{eq:optimal-values-of-alpha}
 \displaystyle  \W^{*} =   \nabla^2_{\W,\W} \JE(\W)^{-1} \cdot \lb  {\basis_0}^T \cdot \B_0^{-1} \cdot \d^{\rm b} + \sum_{k=0}^{\N} \Q_k^T \cdot \R_k^{-1} \cdot \d_{k} \rb\,.
\end{eqnarray}

Since $\x_k^{\as}$ in \eqref{eq:analysis-state-sol-ensemble} represents an approximated solution rather than an exact solution, the initial analysis $\x_0^{\as}$ is only recovered and propagated in time in order to obtain an approximation of the optimal trajectory of \eqref{eq:J4VAR}. 

Equivalent bases for the range of $ \basis_k$ can be utilized to formulate the subspace approximation  \eqref{eq:solution-in-ensemble-space}. For instance, the proper orthogonal decomposition (POD) \cite{Tian2008} is widely used to obtain a basis that captures most of the variance of the snapshot \eqref{eq:snapshots}. Consider the matrix of snapshots deviations 
\begin{eqnarray*}
\displaystyle 
\DXS = \frac{1}{\sqrt{\Nens}}\lb \basis_0\,^T,\basis_1\,^T, \ldots, \basis_{\N}\,^T\rb^T \in \Re^{(\Nstate \cdot (\N+1)) \times \Nens}\,,
\end{eqnarray*}
and its singular value decomposition (SVD)
\begin{eqnarray}
\DXS = {\bf U} \cdot {\boldsymbol \Sigma} \cdot {\bf V} \in \Re^{(\Nstate \cdot (\N+1)) \times \Nens} \,,
\end{eqnarray}
where ${\bf U} \in \Re^{(\Nstate \cdot (\N+1))\times (\Nstate \cdot (\N+1))}$ and ${\bf V} \in \Re^{\Nens \times \Nens}$ are the right and left singular vectors, respectively, and $\Sigma = {\rm diag} \{\sigma_1, \sigma_2,\ldots,\sigma_{\Nens}\} \in \Re^{(\Nstate \cdot (\N+1)) \times \Nens}$ is a diagonal matrix whose diagonal entries are the singular values with $\sigma_1 \ge \sigma_2 \ge \ldots \ge \sigma_{\Nens}$. 
Since
\begin{eqnarray*}
\displaystyle {\DXS}^T \cdot \DXS = {\bf V} \cdot {\bf \Sigma}^2 \cdot {\bf V}^{T} \in \Re^{\Nens \times \Nens} \,,
\end{eqnarray*}
the POD basis vectors can be computed as 
\begin{eqnarray}
\label{eq:POS-basis}
\displaystyle \PODB_k  = \basis_k \cdot {\bf V} \cdot \Sigma^{-1/2} \in \Re^{\Nstate \times \Nens} \,,
\end{eqnarray}
and therefore, equivalent to \eqref{eq:solution-ensemble-space}, $\x_k$ can be expressed as follows:
\begin{eqnarray}
\label{eq:analysis-sol-POD-modes}
\displaystyle \x_k = \xm_k+\sum_{i=1}^{\nPOD} \wPOD_{i} \cdot \lp \frac{\basis_k \cdot {\bf v}_i}{\sqrt{\sigma_i}}\rp = \xm_k+ \PODB_k^{\nPOD} \cdot \WPOD\,,
\end{eqnarray}
where we have chosen the columns of ${\bf \Sigma}$ to be orthonormal, $ \PODB_k^{\nPOD}$ holds the first $\nPOD$ basis vectors, $\WPOD = \lb \wPOD_1,\wPOD_2,\ldots,\wPOD_{\nPOD} \rb^T \in \Re^{\nPOD \times 1}$ is the vector of weights to be determined, and $\nPOD$ can be computed as follows
\begin{eqnarray}
\label{eq:POD-number-of-basis}
\displaystyle \nPOD = \underset{p} {\mathrm{arg\,min}} \, \ll p, I(p):\frac{\sum_{i=1}^p \sigma_i}{\sum_{i=1}^{\Nens} \sigma_i}>\gamma: \gamma \in \lp 0,\,1 \rp \rl \,.
\end{eqnarray}

Note that, the parameter $\gamma$ provides how much variance (sometimes called kinetic energy) we want to retain in the POD bases, commonly the values of $\gamma$ ranges in $(0.9 \,,0.95)$. It is well known that POD bases are the most efficient among all possible linear combinations in the sense, for a given number $\nPOD$ of basis vectors, POD decomposition captures the most possible variance \cite{Ly2001,Ly2002}. In addition, POD bases reduce the equation \eqref{eq:J4VAR-alpha} to
\begin{eqnarray}
\label{eq:POD-function-beta}
\displaystyle \JE^\textnormal{\sc pod}\lp \WPOD \rp = \frac{1}{2} \cdot \Nens \cdot  \ln \WPOD \rn^2  + \frac{1}{2} \cdot \sum_{k=0}^{\N} \ln \d_{k} - \Z_k \cdot \WPOD\rn_{\R_k^{-1}}^2 \,,
\end{eqnarray}
whose first and second derivatives are 
\begin{subequations}
\label{eq:POD-derivatives}
\begin{eqnarray}
\nonumber
\displaystyle \nabla_\WPOD \JE^\textnormal{\sc pod} \lp \WPOD \rp &=& \lb \Nens \cdot \I_{\nPOD \times \nPOD} + \sum_{k=0}^{\N}   \Z_k^T \cdot \R_k^{-1} \cdot \Z_k \rb \cdot \WPOD \\
\label{eq:POD-1st-derivative}
&-& \sum_{k=0}^{\N} \Z_k^T \cdot \R_k^{-1} \cdot \d_{k} \in \Re^{\nPOD \times 1} \,, \\
\label{eq:POD-2st-derivative}
\displaystyle \nabla^2_{\WPOD,\WPOD} \JE^\textnormal{\sc pod} \lp \WPOD \rp &=& \Nens \cdot \I_{\nPOD \times \nPOD} + \sum_{k=0}^{\N}   \Z_k^T \cdot \R_k^{-1} \cdot \Z_k \in \Re^{\nPOD \times \nPOD}\,,
\end{eqnarray}
\end{subequations}
where $\Z_k = \H_k \cdot \PODB_k$ and $\I_{r \times r}$ is the identity matrix of dimension ${r \times r}$. Thus, an equivalent problem to \eqref{eq:solution-ensemble-space} is
\begin{eqnarray}
\label{eq:solution-ensemble-space-POD}
 \WPOD^{*} = \underset{\WPOD} {\mathrm{arg \,min}} \,\JE^\textnormal{\sc pod}(\WPOD) \in \Re^{\nPOD \times 1}\,,
\end{eqnarray}
whose solution reads:
\begin{eqnarray}
\label{eq:optimal-values-of-beta}
\displaystyle 
\WPOD^{*} =  \nabla^2_{\WPOD, \WPOD} \JE^\textnormal{\sc pod} \lp \WPOD \rp^{-1} \cdot \lb \sum_{k=0}^{\N} {\Z^{(k)}}^T \cdot \R_k^{-1} \cdot \d_{k} \rb \,.
\end{eqnarray}

Data assimilation methods that make use of the POD basis (such as, for example, POD-4D-EnKF  \cite{Tian2008}) are defined as follows:
\begin{enumerate}
\item {\bf Ensemble generation}. The ensemble \eqref{eq:initial-ensemble} is built centered at $\x_0^{\rm b}$ with covariance matrix $\B_0$. The ensemble members are propagated and $\N+1$ snapshots of each member are saved.
\item {\bf Basis computation}. The POD basis \eqref{eq:POS-basis} are computed and $\nPOD$ vectors are selected according to \eqref{eq:POD-number-of-basis}.
\item {\bf Compute reduced-space solution}. Compute the optimal weights solution\eqref{eq:optimal-values-of-beta}.
\item {\bf Compute full-space initial condition}. Let $\x_0^{\as} = \x_0^{\bs}+\PODB^{\nPOD}_{0} \cdot \WPOD^{*}$.
\item {\bf Propagate analysis}. $\x_{k}^{\bs} = \M_{t_{k-1} \rightarrow t_{k}}(\x_{k-1}^{\as})$ for $1 \le k \le \N$.
\end{enumerate}

According to Tian \cite{Tian2008} the POD bases capture not only the spatial structure of the state  but also its temporal evolution.

The optimal solution of the POD-4D-EnKF provides an approximation of the analysis \eqref{eq:analysis-full-model}. The process can be continued in an iterative fashion in order to improve the analysis; the solution of one iteration becomes the new background state for the next iteration. The idea of using a sequence of minimizations of the surrogates  \eqref{eq:J4VAR-alpha} or \eqref{eq:POD-function-beta} in order to approach the minimum of \eqref{eq:J4VAR} has been explored in the derivative-free optimization literature \cite{ConnDerivativeFree2009}. 

A rigorous implementation has been recently proposed by Gratton et al. \cite{Gratton2013}. The method is called {\it Iterative Subspace Minimization} (ISM) and solves iteratively the problem \eqref{eq:J4VAR} via the projection of the full space onto the space spanned by the POD bases. The ISM method is defined as follows:
\begin{enumerate}
\item {\bf Initialization}. Let $\x_0^{[0]}  \leftarrow \xm_0$ (the initial background) and $j  \leftarrow 0$. 
\item {\bf Ensemble generation}. The ensemble \eqref{eq:initial-ensemble} is built centered at $\x_0^{[j]}$ with covariance matrix $\B_0$. The ensemble members are propagated and $\N+1$ snapshots are saved.
\item {\bf Basis computation}. The POD basis \eqref{eq:POS-basis} are computed and $\nPOD$ vectors are selected according to \eqref{eq:POD-number-of-basis}.
\item {\bf Suproblem solution}. The optimization problem \eqref{eq:solution-ensemble-space-POD} is partially solved making use of the {\it Coordinate Search Method (CSM)} \cite{Custodio2007,Custodio2010}, from which we obtain $\WPOD^{*}$.
\item {\bf Solution update}. Set ${\x}^{[j+1]}_0 \leftarrow {\x}^{[j]}_0+\Phi_0^{\nPOD} \cdot \WPOD^{*}$, $j \leftarrow j+1$, and go to Step 2.
\end{enumerate}

The ISM method solves the optimization subproblem \eqref{eq:POD-function-beta} via the CSM approach which does not make use of derivative information, and therefore no optimality conditions are checked. Other methods can be used at this step. For instance, one can employ the analytical solution \eqref{eq:solution-ensemble-space-POD}, which guarantees to obtain the local minimizer of each subproblem and reduce the total number of outer iterations and function evaluations. 

The Trust Region (TR) framework can be employed in order to exploit the information brought by the derivatives of the ensemble cost functions \eqref{eq:J4VAR-alpha} and \eqref{eq:POD-function-beta} and to provide descent directions. One of the most attractive features of TR methods is that they are provably globally convergent 
under general assumptions \cite{ConnCH06,ConnCH07,ConnCH08}. A general overview of the TR approach is presented in the appendix \ref{app:TR-method}. To the best of our knowledge TR methods have not been used yet
in the context of ensemble-variational data assimilation. This work develops a TR-based approach which
performs a sequence of optimizations in ensemble spaces. The ensemble based partial solutions are 
linked to the full space solutions at each iteration. The background error statistics of the estimates obtained at each iteration 
are linked to the TR radius size and the spread of the underlying ensemble.
The new method enjoys all these properties and is presented in the next section.

\section{The TR-4D-EnKF method}
\label{sec:proposed-method}
{
In this section we develop a Trust Region 4D-EnKF (TR-4D-EnKF) approach to data assimilation. We start with a general overview of the method and then present the computational algorithm in detail.

The initial solution and background covariance matrix in the model space are given by the initial approximation of the background $\x_0^{[0]} = \xb_0$ and error covariance matrix $\B_0^{[0]} = \B_0$, respectively. The initial ensemble \eqref{eq:initial-ensemble-members} is drawn from $\Nor \lp \x_0^{[0]} \,, \B_0^{[0]} \rp$. The ensemble members are propagated in time and $\N+1$ snapshots are stored as we have discussed previously. We initialize the vector of weights to $\W = {\bf 0}_{\Nens}$ and $\rm j = 0$.  
}
In order to solve the numerical optimization problem \eqref{eq:analysis-full-model}
we build a quadratic model for the cost function \eqref{eq:J4VAR} optimization process. The standard approach makes use of the full space gradient,
and possibly Hessian, of \eqref{eq:J4VAR}. We seek to avoid the implementation of a full adjoint model to compute exact derivatives. The idea is
to approximate the derivatives of $\J(\x)$ by the ensemble space derivatives \eqref{eq:J4VAR-derivat} and Hessian \eqref{eq:J4VAR-derivat2}. The resulting quadratic model is: 
\begin{eqnarray*}
\displaystyle \m^{\rm [j]} \lp \sa \rp &=& \JE(\W+\sa) \\
 &=&\frac{1}{2}\cdot \ln \d^{\rm b} - \basis_0\cdot \lp \W+\sa \rp \rn_{\B_0^{-1}}^2+ \frac{1}{2} \sum_{k=0}^{\N} \ln \d_{k}^{\rm o}- \Q_k \cdot \lp \W+\sa \rp \rn_{\R_k^{-1}}^2\\
&=& \underbrace{\frac{1}{2}\cdot \ln \d^{\rm b} - \basis_0 \cdot \W \rn_{\B_0^{-1}}^2+ \frac{1}{2} \sum_{k=0}^{\N} \ln \d_{k}^{\rm o}- \Q_k \cdot \W \rn_{\R_k^{-1}}^2}_{\JE(\W)} \\
&+& \underbrace{ \ll \lb \basis_0^T \cdot \B_0^{-1} \cdot \basis_0 +  \sum_{k=0}^{\N} \Q_k^T \cdot \R_k^{-1} \cdot \Q_k  \rb \cdot  \W - {\bf c} \rl^T \cdot \sa}_{ \nabla \JE(\W)^T \cdot \sa} \\
&+& \underbrace{\frac{1}{2}  \sa^T  \lb \basis_0^T \cdot \B_0^{-1} \cdot \basis_0 +  \sum_{k=0}^{\N} \Q_k^T \cdot \R_k^{-1} \cdot \Q_k   \rb  \sa}_{\frac{1}{2} \sa^T \nabla^2 \JE(\W) \sa} \,,
\end{eqnarray*}
where ${\bf c} = {\basis_0}^T \cdot \B_0^{-1} \cdot \d^{\bs} +\sum_{k=0}^{\N} \Q_k^T \cdot \R_k^{-1} \cdot \d_{k}  \in \Re^{\Nens \times 1}$. 
This can be rewritten as
\begin{eqnarray}
\label{eq:TR-quadraticmodel-ens}
\displaystyle \m^{\rm [j]} (\sa) = \JE(\W)+\sa^T \nabla_{\W} \JE(\W)+\frac{1}{2} \sa^T \nabla_{\W,\W}^2 \JE(\W) \sa \,.
\end{eqnarray}
The optimal step $\sa^{*}$ in the ensemble space is given by the solution of the constrained optimization sub-problem
\begin{subequations}
\label{eq:subproblem-tr}
\begin{eqnarray}
\label{eq:subproblem-to-solve}
 \sa^{*} &=& \underset{\sa} {\mathrm{arg\,min}} \,\m^{\rm [j]} (\sa)\,, \\
\label{eq:TR-model-space-constraint}
&& \textnormal{subject to  } \displaystyle \|\basis_0 \cdot \lp \W + \sa \rp \| \le \radius^{\rm [j]}  \,.
\end{eqnarray}
\end{subequations}
The trust region constraint \eqref{eq:original-trust-region-constraint} is formulated such as to use the trust region radius $\Delta^{\rm [j]} $ from the full model space. 

The solution of \eqref{eq:subproblem-tr} provides the following trial point in the ensemble space
\begin{subequations}
\begin{eqnarray}
\label{eq:trial-points-ens}
\displaystyle \W^{'} = \W+\sa^{*} \,,
\end{eqnarray}
which corresponds to the following state in the model space
\begin{eqnarray}
\label{eq:trial-points-mod}
\displaystyle \x_0^{\rm '} = \x_0^{[j]} + \underbrace{ \basis_0 \cdot \overbrace{\lp \W+\sa^{*}  \rp}^{\W^{\rm trial}}}_{\dx^{*}} \,.
\end{eqnarray}
\end{subequations}
The problem \eqref{eq:subproblem-tr} is solved using Lagrangian multipliers. The first and second derivatives of the model \eqref{eq:TR-quadraticmodel-ens} are 
\begin{subequations}
\begin{equation}
\label{eq:first-derivative-model}
\displaystyle \nabla\m^{\rm [j]} (\sa) = \nabla \JE(\W)+\nabla^2 \JE(\W) \cdot \sa \in \Re^{\Nens \times 1}\,,
\end{equation}
and
\begin{equation}
\label{eq:second-derivative-model}
\displaystyle \nabla^2\m^{\rm [j]} (\sa) = \nabla^2 \JE(\W) \in \Re^{\Nens \times \Nens}\,,
\end{equation}
\end{subequations}
respectively. The trust region constraint \eqref{eq:TR-model-space-constraint} can be written as 
\begin{eqnarray}
\label{eq:model-space-constrained}
\displaystyle \ln \W + \sa \rn_{\P}^2  - {\radius^{\rm [j]} }^2 + \varsigma^2 =0\,,
\end{eqnarray}
where $\varsigma \in \Re$ is a slack variable and $\P = {\basis_0}^T \cdot {\basis_0} \in \Re^{\Nens \times \Nens} $. Consider the Lagrangian 
\begin{eqnarray}
\label{eq:lagrangian-problem}
\displaystyle \La \lp \sa,\,\lambda,\,\varsigma \rp = \m(\sa)+\lambda \lp \ln \W + \sa \rn_{\P}^2  - {\radius^{\rm [j]} }^2 + \varsigma^2 \rp\,,
\end{eqnarray}
The constrained problem \eqref{eq:subproblem-tr} becomes the unconstrained optimization problem
\begin{eqnarray}
\label{eq:lagrangian-step-size}
\sa^{*} = \underset{\sa} {\mathrm{arg\,min}} \,\La \lp \sa,\,\lambda,\,\varsigma \rp \,.
\end{eqnarray}
The stationarity conditions for \eqref{eq:lagrangian-problem} read:
\begin{eqnarray}
\label{eq:lagrangian-stationary-conditions}
\displaystyle \nabla \La \lp \sa,\,\lambda,\,\varsigma \rp = \begin{bmatrix} \La_{\sa} \lp \sa,\,\lambda,\,\varsigma \rp \\
\La_{\lambda} \lp \sa,\,\lambda,\,\varsigma \rp \\
\La_{\varsigma} \lp \sa,\,\lambda,\,\varsigma \rp  \end{bmatrix} = 0\,,
\end{eqnarray}
where 
\begin{eqnarray*}
\displaystyle 
\La_{\sa} \lp \sa,\,\lambda,\,\varsigma \rp &=& \nabla \m \lp \sa \rp +2 \lambda \,  \P \cdot \lp \W+\sa \rp = \mathbf{0}  \in \Re^{\Nens} \,, \\
 \La_{\lambda} \lp \sa,\,\lambda,\,\varsigma \rp &=&  \ln   \W + \sa \rn_{\P}^2 - {\radius^{\rm [j]} }^2 + \varsigma^2 = 0 \in \Re \,,\\
\La_{\varsigma} \lp \sa,\,\lambda,\,\varsigma \rp &=& 2 \cdot \lambda \cdot \varsigma  = 0 \in \Re\,,
\end{eqnarray*}
which provides all the information needed to solve \eqref{eq:lagrangian-step-size}. Note that, when the full step is taken in the ensemble space
\begin{eqnarray*}
\displaystyle \ln \basis_0  \cdot \W^{\rm '}  \rn \le \Delta^{\rm [j]}  \,,
\end{eqnarray*}
the exact solution \eqref{eq:solution-ensemble-space} can be employed. Then, $\N+1$ snapshots of the full model solution started from $\x_0^{\rm '}$ \eqref{eq:trial-points-mod} are stored. 
The following ratio is computed:
\begin{eqnarray}
\label{eq:rho-trust-region-data-assimilation}
\rho^{\rm [j]}  = \frac{\J \lp \x^{\rm [j]} \rp - \J \lp \x^{\rm '} \rp}{ \m \lp {\bf 0}_{\Nens} \rp - \m \lp \sa^{*} \rp} = \frac{\J \lp \x^{\rm [j]} \rp - \J \lp \x^{\rm '} \rp}{\JE \lp \W \rp - \JE \lp \W^{\rm '}\rp }\,.
\end{eqnarray}
Based on the $\rho^{\rm [j]} $ value, the next updates are made for the solution in the model space
\begin{eqnarray}
\label{eq:solution-update-TR}
 \x^{\rm [j+1]}  &:=& \begin{cases}
 \x^{\rm [j]} & \text{for $\rho^{\rm [j]}  \le \eta$}, \\
 \x^{\rm '}  & \text{otherwise}, 
\end{cases} 
\end{eqnarray} 
and for the TR radius size
\begin{eqnarray}
\label{eq:radius-update-TR}
\radius^{\rm [j+1]} &:=& \begin{cases}
\radius^{\rm [j]} \cdot \gamma_{\rm dec}  & \text{for $\rho^{\rm [j]}  < \theta_1$}, \\
\radius^{\rm [j]} & \text{for $\theta_1 \le \rho^{\rm [j]}  < \theta_2$ or $\rho^{\rm [j]} >1$}, \\
\min{\lp \radius^{\rm [j]} \cdot \gamma_{\rm inc},\,\radius_{\rm max} \rp} & \text{for $\theta_2 \le \rho^{\rm [j]}  \le 1$}.
\end{cases} 
\end{eqnarray}

The current solution becomes the new background and therefore, a new ensemble of full model solutions is generated, snapshots are taken, a new set of basis vectors is built, and the overall process is repeated.

Since a partial assimilation of observations has been carried out the uncertainty associated with the new background is changed. As an analogy, in the EnKF the spread of the ensemble members around the background is decreased after the analysis step. Consequently, before generating a new ensemble, we want to adjust the spread of the background errors. This is done according to the heuristic formula
\begin{eqnarray}
\label{eq:tr-B-value}
\displaystyle \B_0^{\rm [j+1]} := \lambda_{\B}(\radius) \cdot \B_0^{\rm [j]}\,,
\end{eqnarray}
where $\lambda_{\B} \lp \Delta \rp$ is a function of the current TR radius size. Note that the TR radius is large when the decrease of the current (quadratic) model is a good predictor of the full model function decrease. In our context, if the dynamics of the full (nonlinear) model is well represented by the ensemble, the prediction done using the quadratic model $\m(\sa)$ is close to the actual reduction of the cost function $\J\lp \x \rp$ and the TR radius is increased. In this case, we want the $\lambda_{\B} \lp \Delta \rp$ value to be small in order to decrease the uncertainty  of the new ensemble around $\x_0^{\rm b}$. Vice-versa, a small TR radius indicates that the current set of basis vectors does not represent well the dynamics of the model. The current assimilation step is not expected to decrease uncertainty; to keep the same uncertainty level for the next ensemble generation we need $\lambda_{\B} \lp \Delta \rp \approx 1$. Both cases are captured by the following heuristic function
\begin{eqnarray}
\label{eq:TR-shrink-value}
\displaystyle 
\lambda_{\B}\lp \radius \rp  = \frac{\radius_{\rm max}}{\radius_{\rm max}+\radius} \,,
\end{eqnarray}
which provides an inverse relation between the TR radius and the spread of the ensemble members. Other functions can be considered as well. In summary,  when the TR radius is large the confidence in the current solution is increased
\begin{eqnarray*}
\lim_{\radius \rightarrow \radius_{\rm max}}
\displaystyle \lambda_{\B} = \frac{\radius_{\rm max}}{2\cdot \radius_{\rm max}} = \frac{1}{2} \,.
\end{eqnarray*}
 On the other hand, when the TR size is small, the current level of background uncertainty remains unchanged for the new ensemble generation
\begin{eqnarray*}
\lim_{\radius \rightarrow 0}
\displaystyle \lambda_{\B} = \frac{\radius_{\rm max}}{\radius_{\rm max}} = 1 \,.
\end{eqnarray*}

\begin{figure}[H]
\centering
\begin{subfigure}{0.5\textwidth}
  \centering
  \includegraphics[width=\linewidth,height=\linewidth]{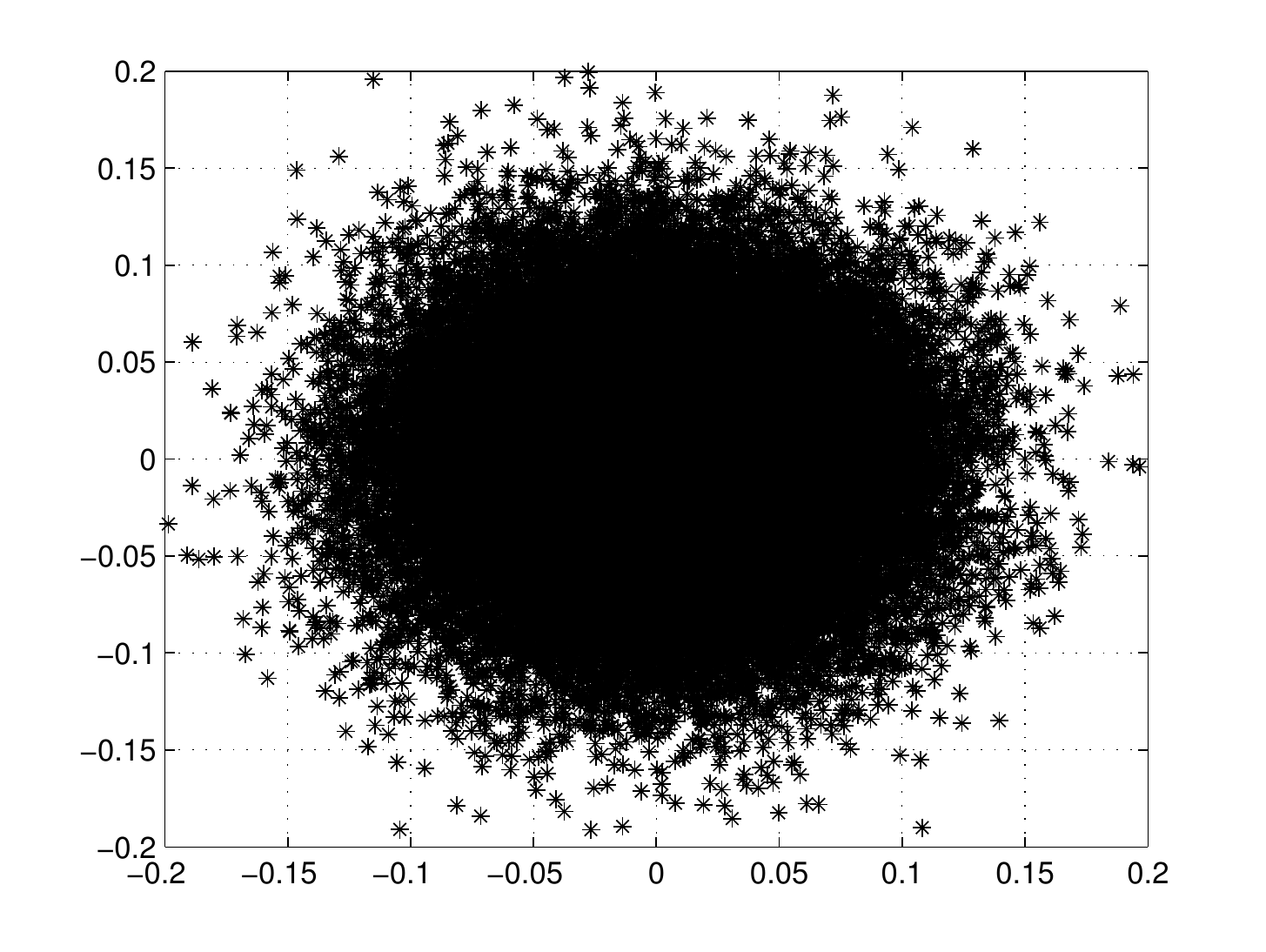}
  \caption{$\lambda_{\B} = 1$}
  \label{fig:B1-lambda-1}
\end{subfigure}%
\begin{subfigure}{0.5\textwidth}
  \centering
  \includegraphics[width=\linewidth,height=\linewidth]{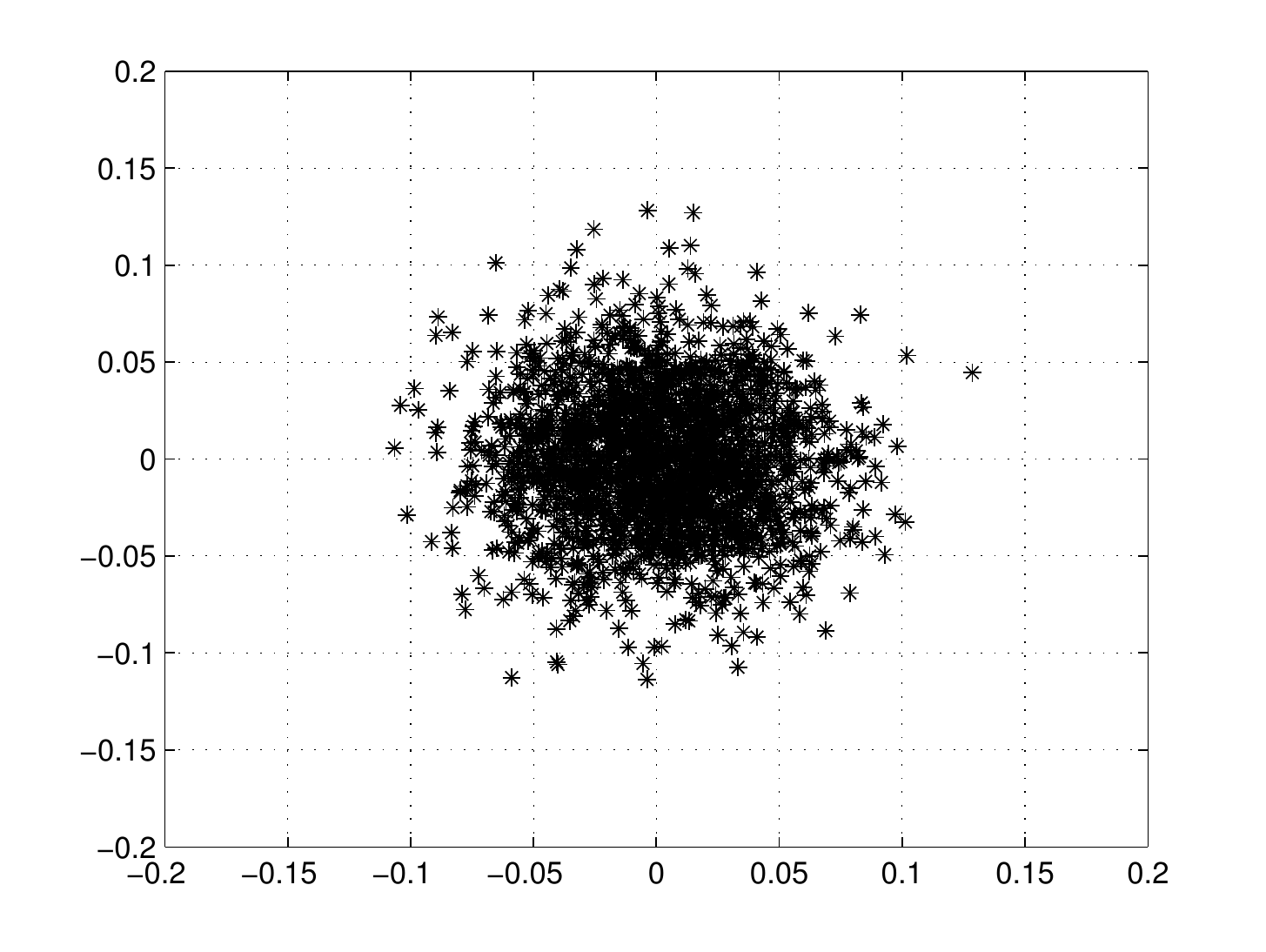}
  \caption{$\lambda_{\B} = \frac{1}{2}$}
  \label{fig:B2-lambda-1_2}
\end{subfigure}
\caption{Impact of the scaling of $\B_0$ on the spread of the newly generated ensemble.}
\label{fig:impact-lambda}
\end{figure}

The effects of the scaling of $\B_0$ on the new background ensemble are shown in the Figure \ref{fig:impact-lambda} for a 2D example. 
The choice $\lambda_{\B}=1$ keeps the uncertainty unchanged (Figure \ref{fig:B1-lambda-1}),
while $\lambda_{\B}=1/2$ shrinks the spread by half (Figure \ref{fig:B2-lambda-1_2}).

The outline of the TR-4D-EnKF method is shown below.{
\begin{enumerate}
\item {\bf Initialization.} The TR parameters are initialized. $\x_0^{\rm [0]} := \x_0^{\rm b}$, $\B_0^{[0]} = \B_0$, $\W = {\bf 0}_{\Nens}$ and $\rm j = 0$.
\item {\bf Ensemble generation.} The initial ensemble is drawn from a distribution with mean $\x_0^{\rm [j]}$ and background error covariance matrix $\B_0^{\rm [j]} $. $\N$ snapshots for each ensemble member are stored \eqref{eq:snapshots} and the basis functions \eqref{eq:matrix-of-anomalies-k} are computed.
\item {\bf Model construction.} Build the quadratic model \eqref{eq:TR-quadraticmodel-ens}.
\item {\bf Optimization problem.} Solve the optimization sub-problem \eqref{eq:subproblem-tr} and compute
$\W^{\rm '}$ and $\x_0^{\rm '}$.
\item {\bf Ratio of prediction.} Run the full model to obtain $\N$ snapshots of the solution 
started from $\x_0^{\rm '}$ and compute the ratio $\rho^{\rm [j]} $ \eqref{eq:rho-trust-region-data-assimilation}.
\item {\bf Solution and TR size update.} Update the solution in the model space and the 
TR radius according to \eqref{eq:solution-update-TR} and \eqref{eq:radius-update-TR}, respectively. 
\item {\bf Background update.} Set $\W := {\bf 0}_{\Nens}$, scale the covariance matrix \eqref{eq:tr-B-value}, $\rm j \gets j+1$, and go to Step 2.
\end{enumerate}

The iteration stopping criterion for the 4D-TR-EnKF implementation can be based on the total number of iterations, on the trace of $\B_0^{\rm [j]}$, or on the trust region radius $\radius^{\rm [j]}$. }

Now we are ready to test our implementation and compare it with other 4D-EnKF implementations discussed in Section \ref{sec:preliminaries}.

\section{Numerical experiments}
\label{sec:experimental}
%
{
In this section we study the accuracy and performance of the TR-4D-EnKF approach. The proposed implementation is compared with the 4D-EnKF implementations discussed in section \ref{sec:preliminaries}: POD and ISM, using the Atmospheric General Circulation Model (AGCM), better known as SPEEDY \cite{SPEEDY,SPEEDY_2,SPEEDY_3,SPEEDY_4}, by the International Centre for Theoretical Physics (ICTP) in Trieste, Italy. The ICTP AGCM  is based on a spectral dynamical core developed at the Geophysical Fluid Dynamics Laboratory (GFDL) \cite{Suarez_SPEEDY} at the Princeton Univeristy Forrestal Campus in Princeton, USA. All the physics are developed on the sphere concisely, the Earth. It is a hydrostatic, $\sigma$-coordinate, spectral-transform model in the vorticity-divergence form described by Bourke \cite{Bourke_SPEEDY}, with semi-implicit treatment of gravity waves. The basic prognostic variables are shown in the Table \ref{tab:instances-used}
\begin{table}[H] \centering
\begin{tabular}{|l|c|c|}\hline
{\bf Variable} & {\bf Symbol} & {\bf Number of Layers} \\\hline
Vorticity & $\Omega$ & 8 \\\hline
Divergence & $\chi$ & 8  \\\hline
Pressure & $\psi$ & 8 \\\hline
Specific Humidity & $\Lambda$ & 8  \\\hline
Temperature & $\tau$ & 1 \\\hline
\end{tabular}
\caption{Pronostic variables in the SPEEDY model.}
\label{tab:instances-used}  
\end{table}

The number of longitudinal and latitudinal points are 96 and 48, respectively. The longitudinal values ranges evenly in $[0,\,96]$ while the latitudinal values ranges evenly in $[0,\,48]$.  This provides a total number of 4096 points per layer. Each layer is mapped to the vector state which provides a total number of $\Nstate = 152064$ components. Only 50\% of the components are observed at each layer, this corresponds to 2048 components per layer and a total number of $\Nobs = 67584$ components being observed in the system.

The metrics used in the tests are the CPU time (which is reported per iteration) and the root mean square error 
\begin{eqnarray}
\label{sec:rmse}
\displaystyle \text{RMSE} = \sqrt{ \frac{1}{\N} \cdot \sum_{k=0}^{\N}  \lp \xt_k-\x^{\as}_k \rp^T \cdot  \lp \xt_k-\x^{\as}_k \rp} \,,
\end{eqnarray}
which provides the average of the squared root differences between the reference solution $\xt$ and the analysis $\x^{\as}$ over the observation times. 

Some details regarding the numerical implementation of the data assimilation methods:
\begin{itemize}
\item Three computational languages are used to carried out the different steps of the compared methods: C, FORTRAN and MATLAB. 
\item The forecast step of the ensemble members is performed in C language making use of MPI. Each ensemble member is independently propagated in time. 
\item The number of ensemble members matches the number of processors.
\item The initial conditions of the ensemble are coded in FORTRAN 77.
\item The assimilation step is carried out in MATLAB.
\item The communication between different programming languages is performed in the data level via NetCDF files.
\item The main core of the program is written in bash language which integrates the different components of the data assimilation process.
\item A two day assimilation window is set-up with observations taken evenly each 3 hours.
\end{itemize}

Other parameters of the numerical simulation are described below.
\begin{itemize}
\item Starting in rest, the ICTP AGMC model is propagated in time for three months, after that, we assume the final state to be the true initial condition $\x^{\rm true}_0$ for our testing
\item Four linear observation operators on the Earth are considered, they are evenly and sequentially distributed over the assimilation window. This mimics, for instance, the use of different sets of sensors at different times in the ocean. The observational operators are shown in figure \ref{fig:linear-operators}.
\begin{figure}[H]
\centering
\begin{subfigure}{.5\textwidth}
  \centering
  \includegraphics[width=1\linewidth,height=1\linewidth]{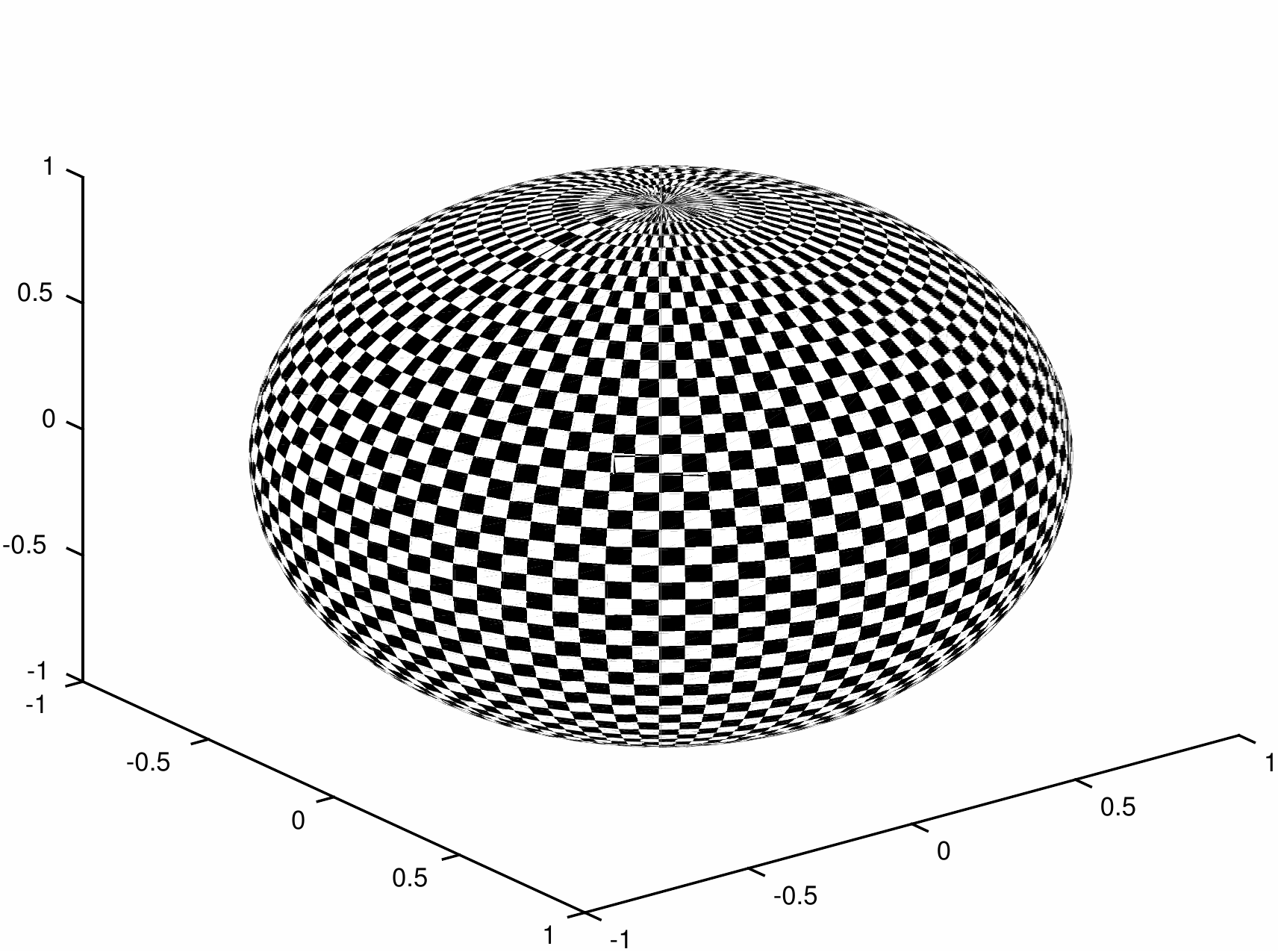}
  \caption{$\H^1$}
  \label{fig:H1}
\end{subfigure}%
\begin{subfigure}{.5\textwidth}
  \centering
  \includegraphics[width=1\linewidth,height=1\linewidth]{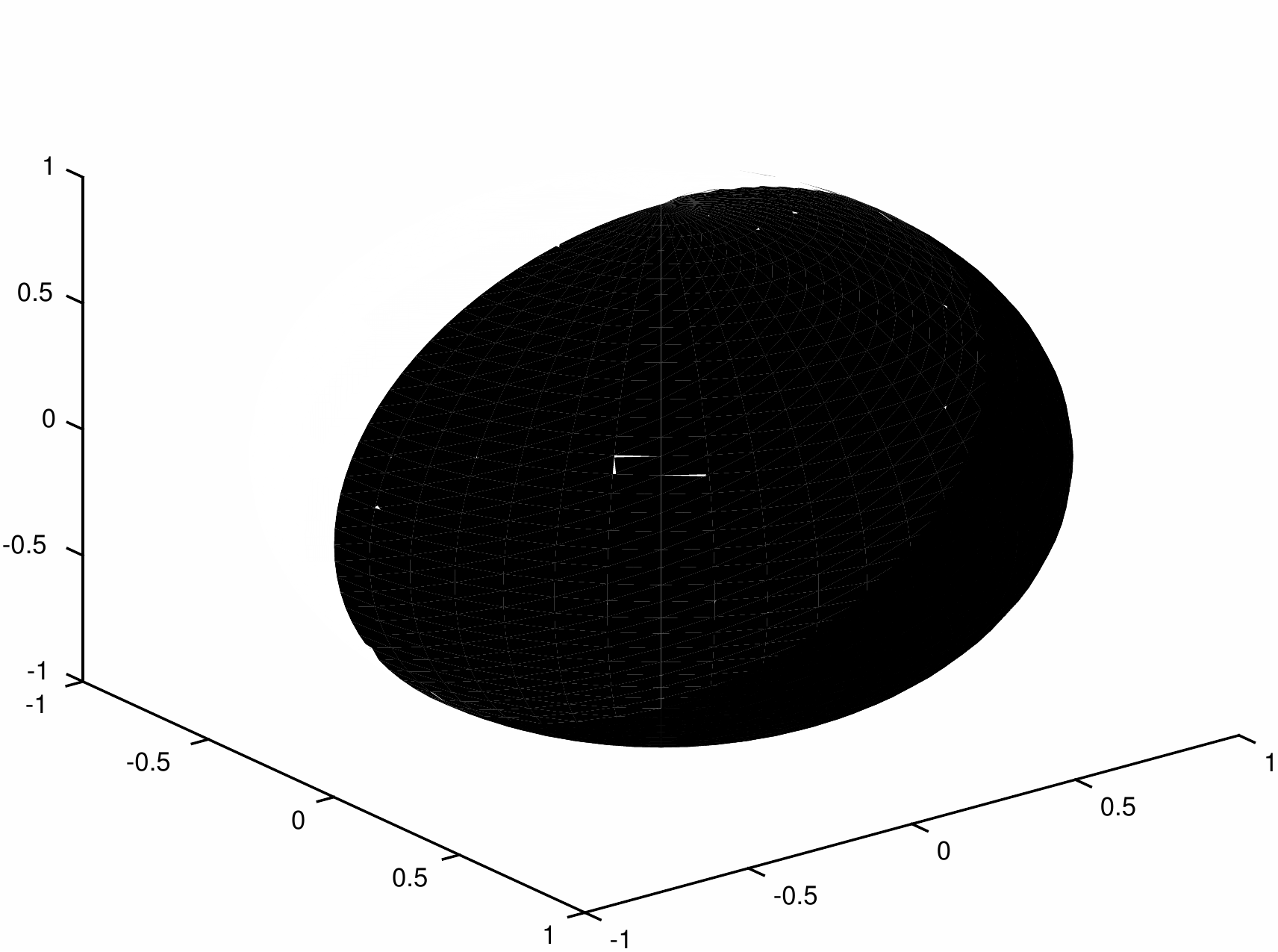}
  \caption{$\H^2$}
  \label{fig:H2}
\end{subfigure}
\begin{subfigure}{.5\textwidth}
  \centering
  \includegraphics[width=1\linewidth,height=1\linewidth]{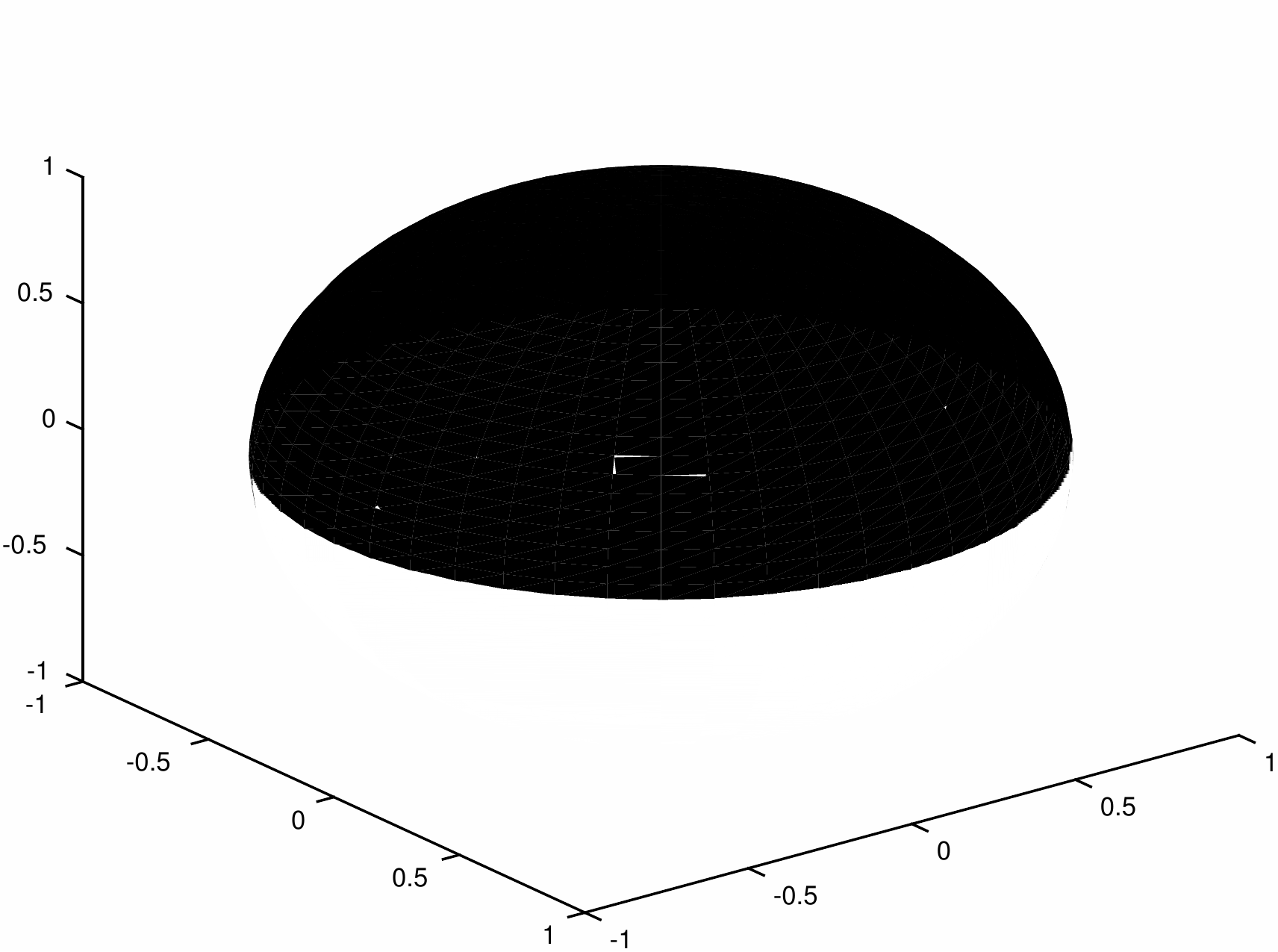}
   \caption{$\H^3$}
  \label{fig:H3}
\end{subfigure}%
\begin{subfigure}{.5\textwidth}
  \centering
  \includegraphics[width=1\linewidth,height=1\linewidth]{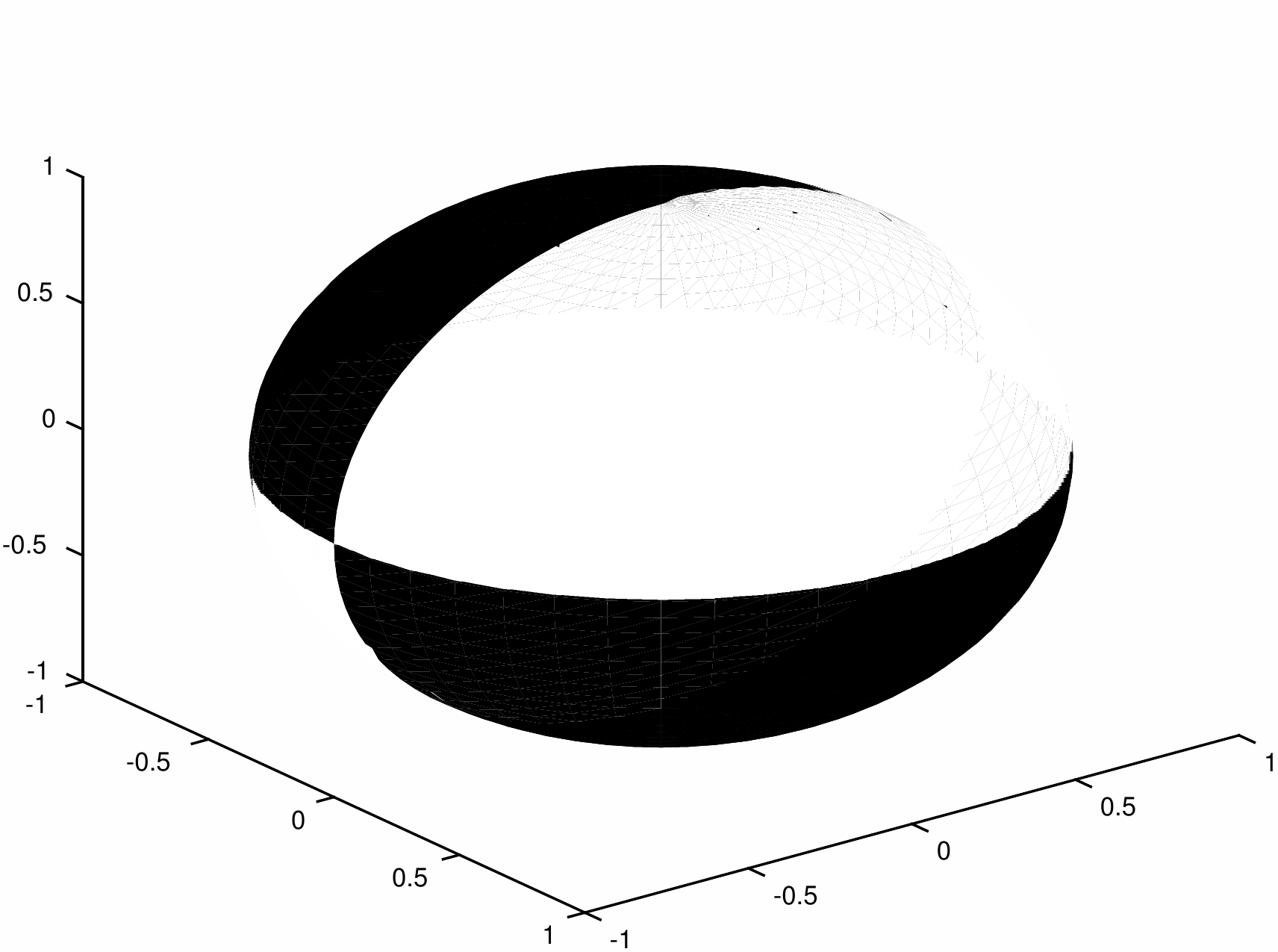}
   \caption{$\H^4$}
  \label{fig:H4}
\end{subfigure}
\caption{Linear observation operators. The dark areas correspond to observed components.}
\label{fig:linear-operators}
\end{figure}
\item The initial background state is a sample from the distribution:
\begin{eqnarray*}
\xb \sim \Nor \lp \x^{\rm true},\, \B_0 \rp \,,
\end{eqnarray*}
where $\B_0 = (0.05)^2 \cdot {\bf I}_{\Nstate \times \Nstate}$. Similarly, the initial ensemble are samples drawn from the distribution
\begin{eqnarray*}
\x^{b(i)} \sim \Nor \lp \xb,\, \B_0 \rp \,,
\end{eqnarray*}
for $1 \le i \le \Nens$.
\item Four ensemble sizes $\Nens$ are considered: 10, 20, 40 and 80.
\item Data errors are normally distributed with parameters
\begin{eqnarray*}
\errobs \sim \Nor \lp {\bf 0}_{\Nobs},\, \R \rp
\end{eqnarray*}
with $\R = (0.01)^2 \cdot {\bf I}_{\Nobs \times \Nobs}$.
\item Five iterations are considered for the ISM and TR-4D-EnKF implementations. This is more than reasonable, in practice, we are not able to propagate the model many times.
\item The parameters for the TR-4D-EnKF optimization are $\gamma_{\rm inc}=1.4$, $\gamma_{\rm dec} = 0.5$, $\radius_{\rm max}=100$, $\radius_0 = 0.1$, $\eta = 0.1$, $\theta_1 = 0.25$ and $\theta_2 = 0.75$.
\end{itemize}

The RMSE and the elapsed times are presented in the Tables \ref{tab:root-mean-square-errors} and \ref{tab:assimilation-times}, respectively. Notice, all the 4D-EnKF methods presented in this paper are able to improve the background initial condition $\x_0^b$ in terms of RMSE. This behaviour holds even in difficult scenarios such as small ensemble sizes (i.e., $\Nens=10$), high dimensional vector states and only 50\% of observed components from the dynamical system. Notice, the POD-4D-EnKF outperforms, in average, the accuracy of the background state by one order of magnitude. The improving is more notorious when the dimension of the ensemble is increased. This obeys to a better representation of the background error statistics onto the space spanned by the ensemble perturbations. As is expected, the more ensemble members, the lesser sampling errors are involved in the assimilation step. Likewise, the ISM method provides very good results with just five iterations and within a reasonable computational effort. The results are much better than the ones obtained via POD-4D-EnKF (equivalent to a single iteration of the ISM). This justifies the iterative refinement of solutions in the context of reduce space approaches. The proposed TR-4D-EnKF outperforms the initial solution in the context of RMSE and after five iterations, the RMSE values look similar to those obtained via the ISM. The figures \ref{fig:snapshots-initial-conditions-vorticity} and \ref{fig:snapshots-initial-conditions-temperature} show the estimated state obtained by each method for the vorticity and the temperature at the Earth's surface. There is no doubt that the initial background state $\x^b_0$ provides a poor estimation of the vorticity and temperature at the Earth's surface.  On the contrary, the solutions obtained by the POD-4D-EnKF seems to be relatively ``close'' to the true state of the system. However, the ISM provides a much better approximation than the POD-4D-EnKF. Likewise, the TR-4D-EnKF approximation is similar to the ISM one. Figure \ref{fig:vorticity-per-ensemble-size-time} provides the time evolution of the errors for the vorticity field among the 8 layers. Note that, the RMSE values in time of the ISM and TR-4D-EnKF are relatively close under the RMSE metric.  However, figure \ref{fig:vorticity-per-ensemble-size} shows a different perspective of this two methods. We report the background and POD-4D-EnKF solutions for comparison purposes since those methods does not require iterations. Note that, in all the cases, the TR-4D-EnKF method performs better than the ISM in the first three iterations. The performance is much better when large ensemble sizes are used in the assimilation window. Note that, the accuracy obtained by five iterations of the ISM is equivalent to that obtained by two iterations of the TR-4D-EnKF. In practice, this is extremely important since model propagation is a labor-intensive process and therefore, the lesser number of times the model is propagated, the better. This implies we are able to obtain good approximations with lesser number of iterations making use of the TR-4D-EnKF than the ISM. This gap between the two methods can be explained in terms of the initial background distribution at each iteration: the TR-4D-EnKF decreases the uncertainty according to trust region sizes, when a good representation of the background error statistics (ensemble members) is contained in the ensemble, the next iteration of the TR-4D-EnKF method will solve an optimization problem where the uncertainty around the initial condition has been decreased. This is not the case of the ISM method where the same uncertainty is hold among all the iterations.

\begin{table}
\begin{center}
    \begin{tabular}{ | c | c | c | c | c | c | c | }
    \hline
    ${\bf \Nens}$ & {\bf Method} & $\boldsymbol{\Omega}$ & $\boldsymbol{\chi}$ & $\boldsymbol{\psi}$ & $\boldsymbol{\Lambda}$ & $\boldsymbol{\tau}$ \\ \hline
    N/A & Background & $1.64(-3)$ & $1.62(-3)$ & $5.58(2)$ & $2.08(3)$ & $1.05(2)$ \\ \hline
    \multirow{3}{*}{10} & POD-4D-EnKF & $9.89(-4)$ & $9.79(-4)$ & $3.29(2)$ & $1.23(3)$ & $6.33(1)$  \\ \cline{2-7}
    & ISM & $5.70(-4)$ & $5.56(-4)$ & $1.82(2)$ & $6.83(2)$ & $3.70(1)$  \\ \cline{2-7}
    & TR-4D-EnKF & $5.66(-4)$ & $5.53(-4)$ & $1.81(2)$ & $6.79(2)$ & $3.66(1)$  \\ \hline
    \multirow{3}{*}{20} & POD-4D-EnKF & $9.07(-4)$ & $8.91(-4)$ & $2.97(2)$ & $1.11(3)$ & $5.92(1)$  \\ \cline{2-7}
    & ISM & $4.07(-4)$ & $4.02(-4)$ & $1.29(2)$ & $4.88(2)$ & $2.67(1)$  \\ \cline{2-7}
    & TR-4D-EnKF & $4.07(-4)$ & $4.02(-4)$ & $1.29(2)$ & $4.87(2)$ & $2.67(1)$  \\ \hline
    \multirow{3}{*}{40} & POD-4D-EnKF & $8.88(-4)$ & $8.80(-4)$ & $2.88(2)$ & $1.07(3)$ & $5.65(1)$  \\ \cline{2-7}
    & ISM & $3.03(-4)$ & $3.10(-4)$ & $9.79(2)$ & $3.69(2)$ & $2.07(1)$  \\ \cline{2-7}
    & TR-4D-EnKF & $2.93(-4)$ & $3.02(-4)$ & $9.47(2)$ & $3.56(2)$ & $2.01(1)$  \\ \hline
    \multirow{3}{*}{80} & POD-4D-EnKF & $7.99(-4)$ & $8.14(-4)$ & $2.67(2)$ & $9.98(2)$ & $5.11(1)$  \\ \cline{2-7}
    & ISM & $2.31(-4)$ & $2.43(-4)$ & $7.49(2)$ & $2.83(2)$ & $1.59(1)$  \\ \cline{2-7}
    & TR-4D-EnKF & $2.11(-4)$ & $2.29(-4)$ & $6.99(2)$ & $2.65(2)$ & $1.46(1)$  \\ \hline
    \end{tabular}
\end{center}
\caption{Root Mean Square Error for different ensemble sizes and data assimilation methods. The notation reads $x(y) = x \times 10^{y}$.}
\label{tab:root-mean-square-errors}
\end{table}

\begin{table}
\begin{center}
    \begin{tabular}{ | c | c | c |}
    \hline
    ${\bf \Nens}$ & {\bf Method} & {\bf Assimilation Time} \\ \hline
    \multirow{3}{*}{10} & POD-4D-EnKF & $\sim$ 10 seconds  \\ \cline{2-3}
    & ISM & $\sim$ 55 seconds  \\ \cline{2-3}
    & TR-4D-EnKF & $\sim$ 1 minute  \\ \hline
    \multirow{3}{*}{20} & POD-4D-EnKF & $\sim$ 12 seconds  \\ \cline{2-3}
    & ISM & $\sim$ 1 minute   \\ \cline{2-3}
    & TR-4D-EnKF & $\sim$ 1 minute  \\ \hline
    \multirow{3}{*}{40} & POD-4D-EnKF & $\sim$ 20 seconds  \\ \cline{2-3}
    & ISM & $\sim$ 1.5 minutes   \\ \cline{2-3}
    & TR-4D-EnKF & $\sim$ 1.9 minutes  \\ \hline
    \multirow{3}{*}{80} & POD-4D-EnKF & $\sim$ 1 minute  \\ \cline{2-3}
    & ISM & $\sim$ 4 minutes  \\ \cline{2-3}
    & TR-4D-EnKF & $\sim$ 5 minutes  \\ \hline
    \end{tabular}
\end{center}
\caption{Assimilation times for the compared 4D-EnKF implementations.}
\label{tab:assimilation-times}
\end{table}
}
\begin{figure}[H]
\centering
\begin{subfigure}{.5\textwidth}
  \centering
    \includegraphics[width=0.7\linewidth,height=0.7\linewidth]{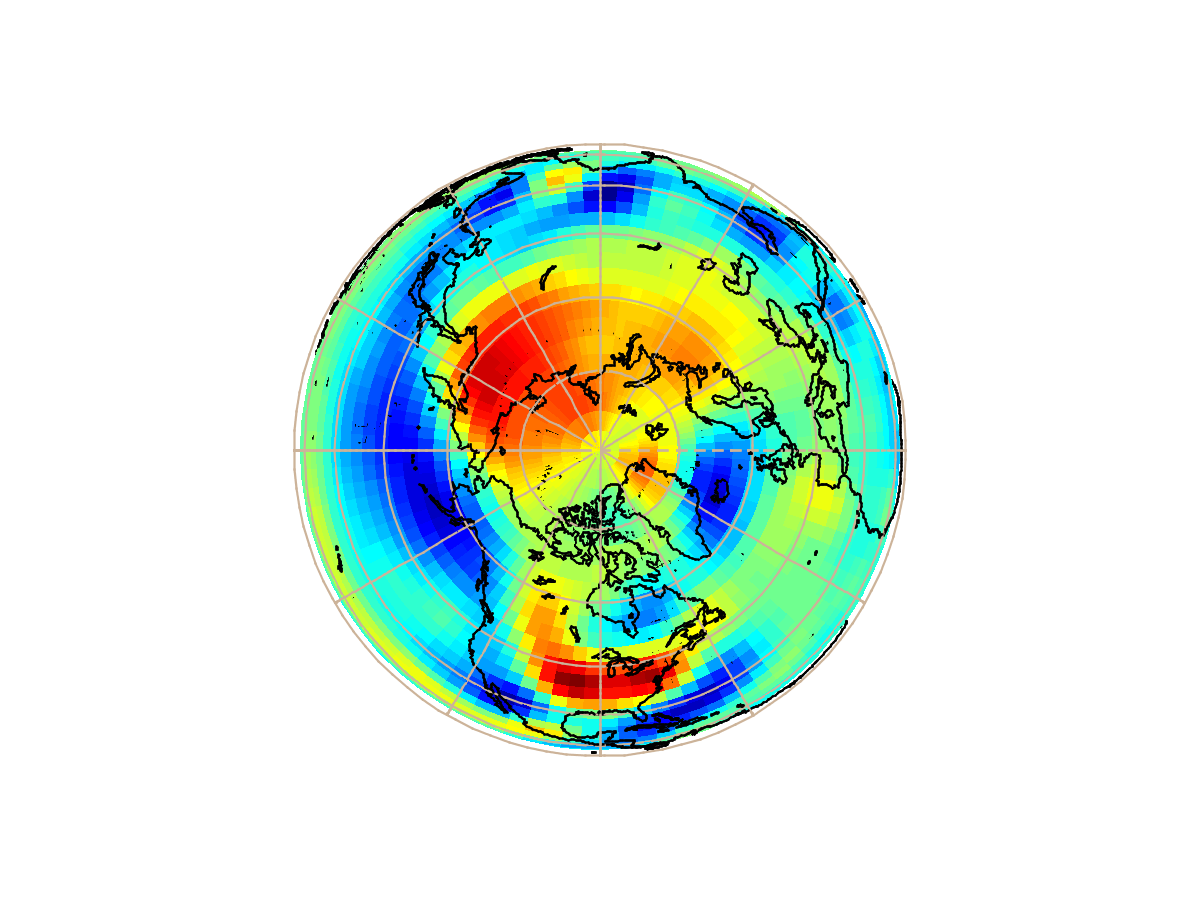}
  \caption{True State}
  \label{fig:V0}
\end{subfigure}\hspace*{-7.8em}%
\begin{subfigure}{.5\textwidth}
  \centering
  \includegraphics[width=0.7\linewidth,height=0.7\linewidth]{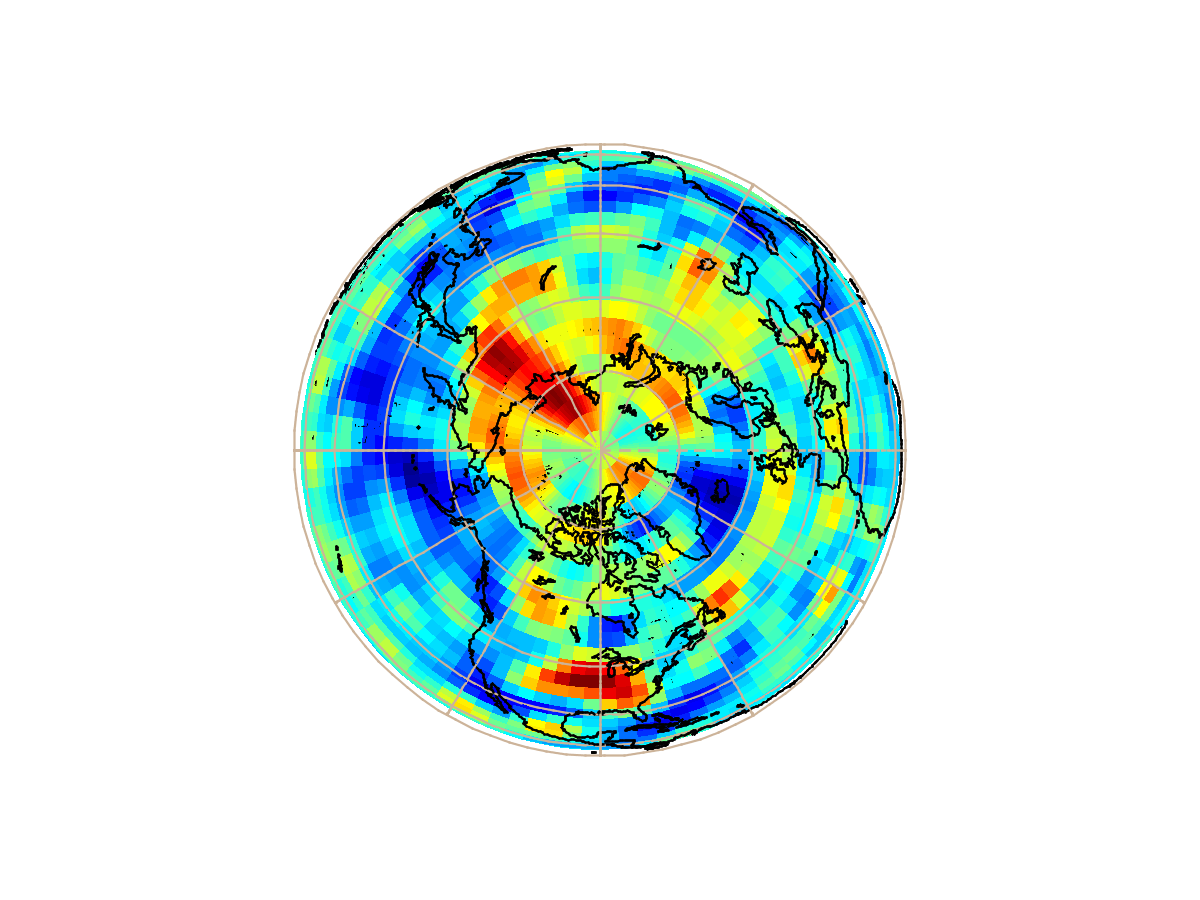}
  \caption{Background}
  \label{fig:V1}
\end{subfigure}\hspace*{-7.8em}%
\begin{subfigure}{.5\textwidth}
  \centering
  \includegraphics[width=0.7\linewidth,height=0.7\linewidth]{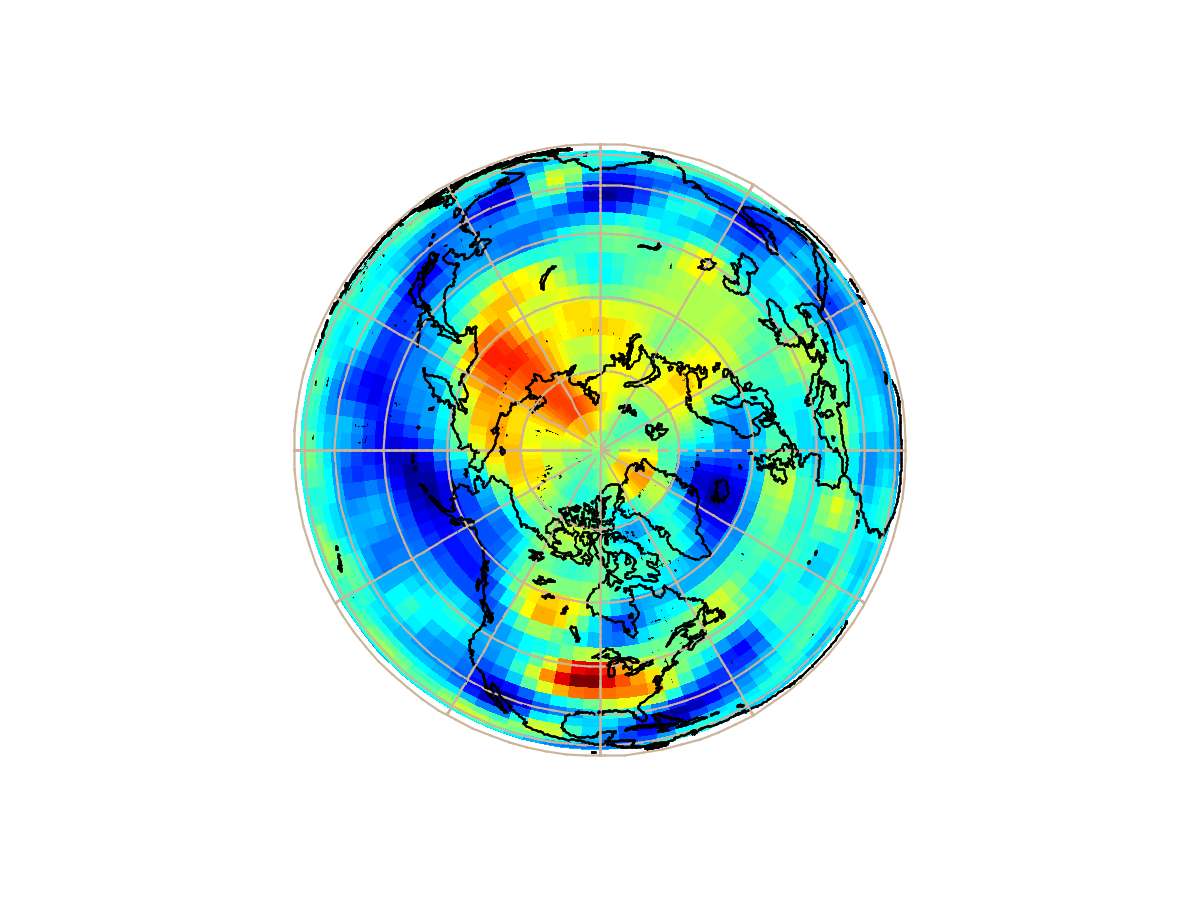}
  \caption{POD-4D-EnKF}
  \label{fig:V2}
\end{subfigure}
\begin{subfigure}{.5\textwidth}
  \centering
  \includegraphics[width=0.7\linewidth,height=0.7\linewidth]{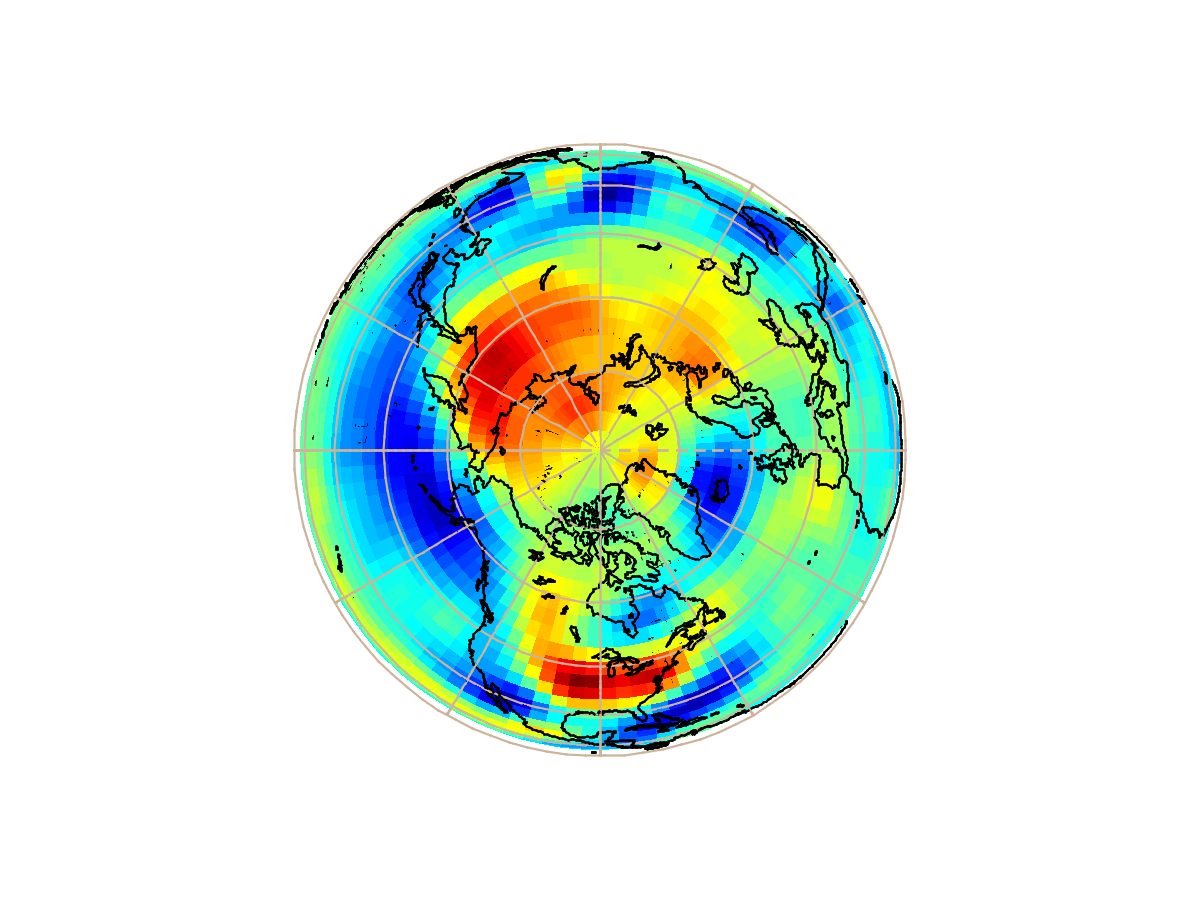}
   \caption{ISM}
  \label{fig:V3}
\end{subfigure}\hspace*{-7.8em}%
\begin{subfigure}{.5\textwidth}
  \centering
  \includegraphics[width=0.7\linewidth,height=0.7\linewidth]{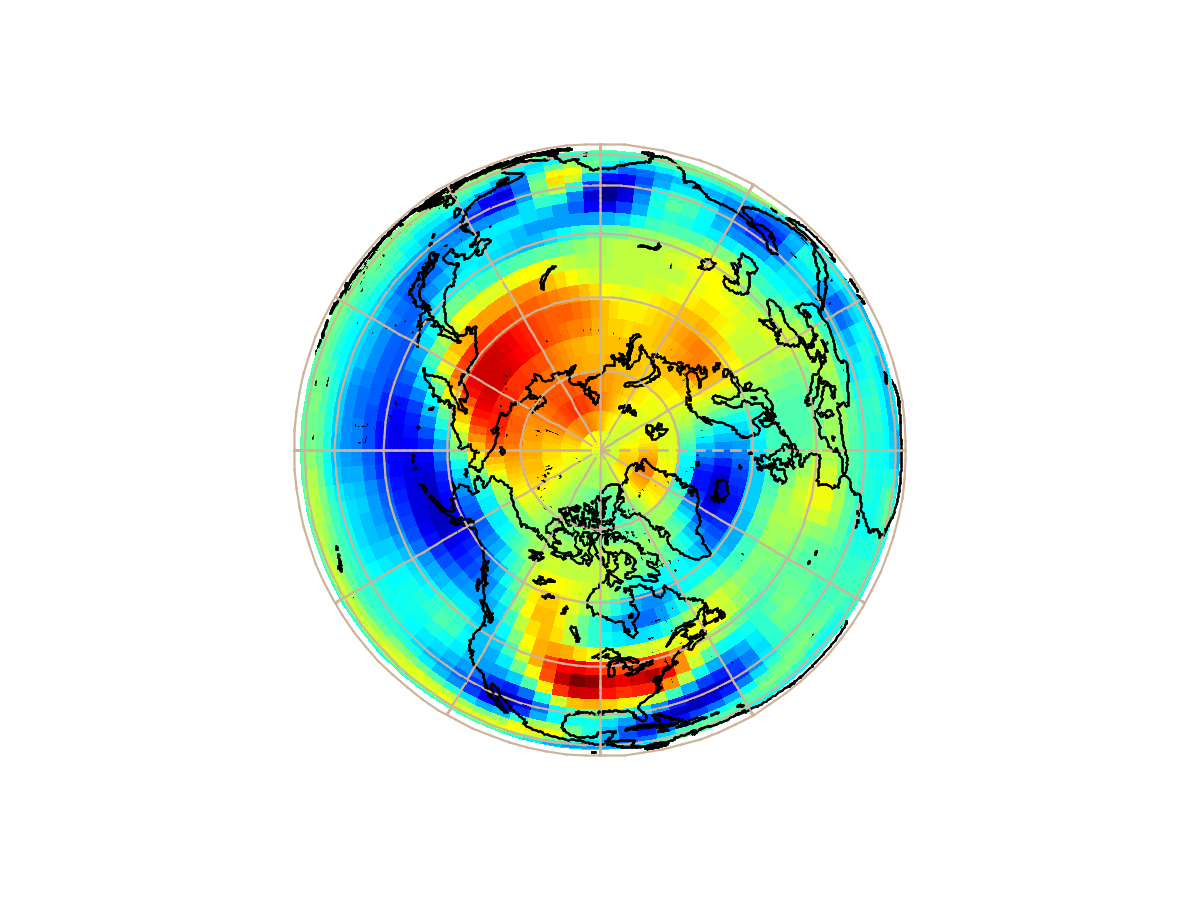}
   \caption{TR-4D-EnKF}
  \label{fig:V4}
\end{subfigure}
\caption{Initial vorticities at the Earth's surface from the analysis states $\x^{a}_0$ for the different compared data assimilation methods. After five iterations, the iterative methods ISM and TR-4D-EnKF provide the most accurate results among the compared implementations.}
\label{fig:snapshots-initial-conditions-vorticity}
\end{figure}

\begin{figure}[H]
\centering
\begin{subfigure}{.5\textwidth}
  \centering
    \includegraphics[width=0.7\linewidth,height=0.7\linewidth]{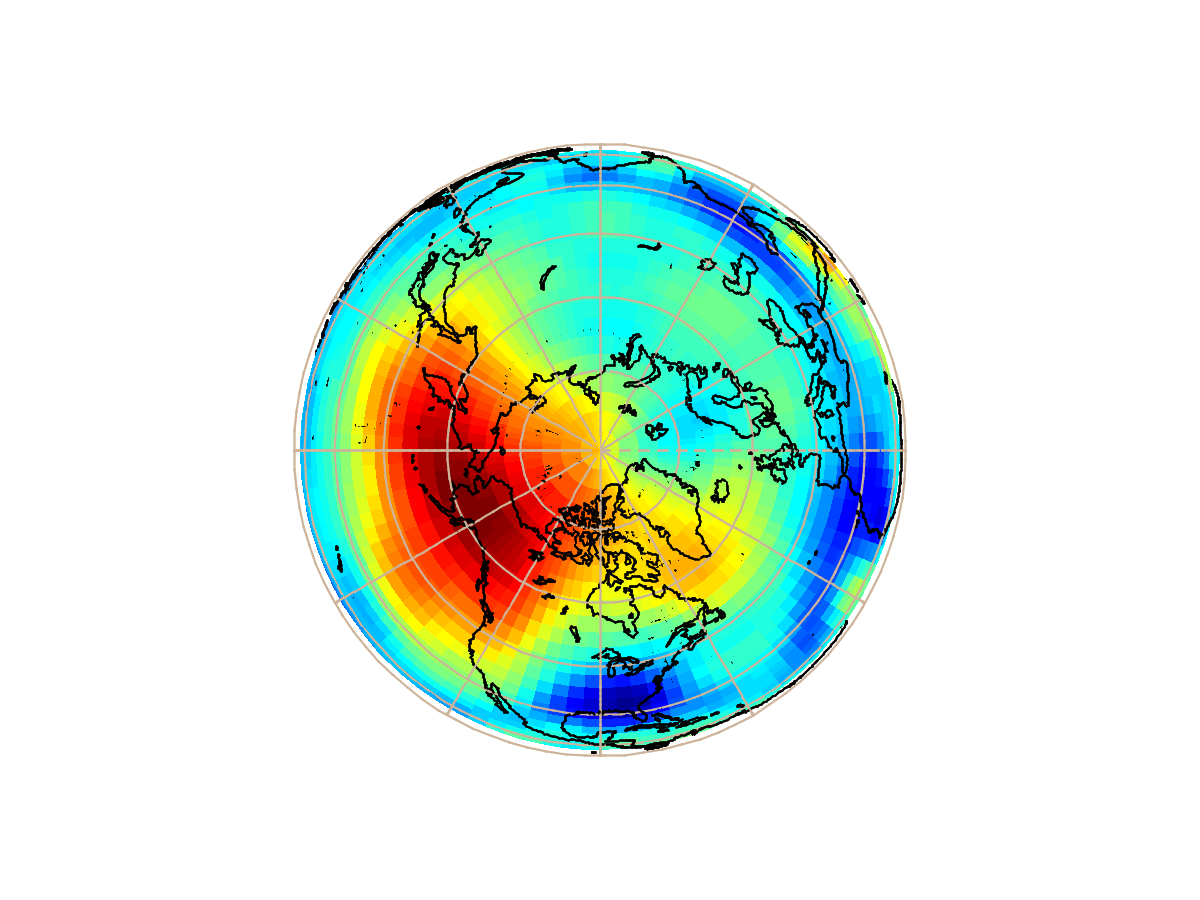}
  \caption{True State}
  \label{fig:T0}
\end{subfigure}\hspace*{-7.8em}%
\begin{subfigure}{.5\textwidth}
  \centering
  \includegraphics[width=0.7\linewidth,height=0.7\linewidth]{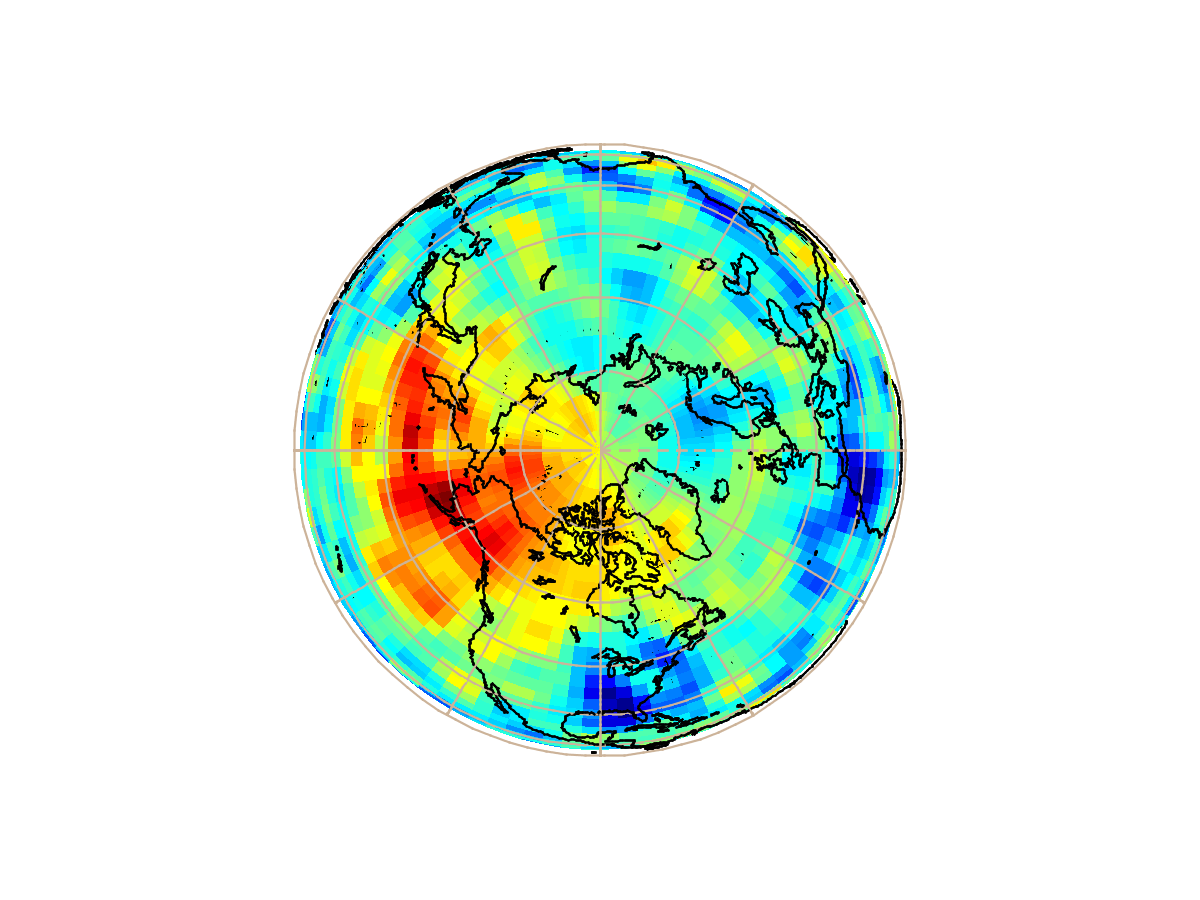}
  \caption{Background}
  \label{fig:T1}
\end{subfigure}\hspace*{-7.8em}%
\begin{subfigure}{.5\textwidth}
  \centering
  \includegraphics[width=0.7\linewidth,height=0.7\linewidth]{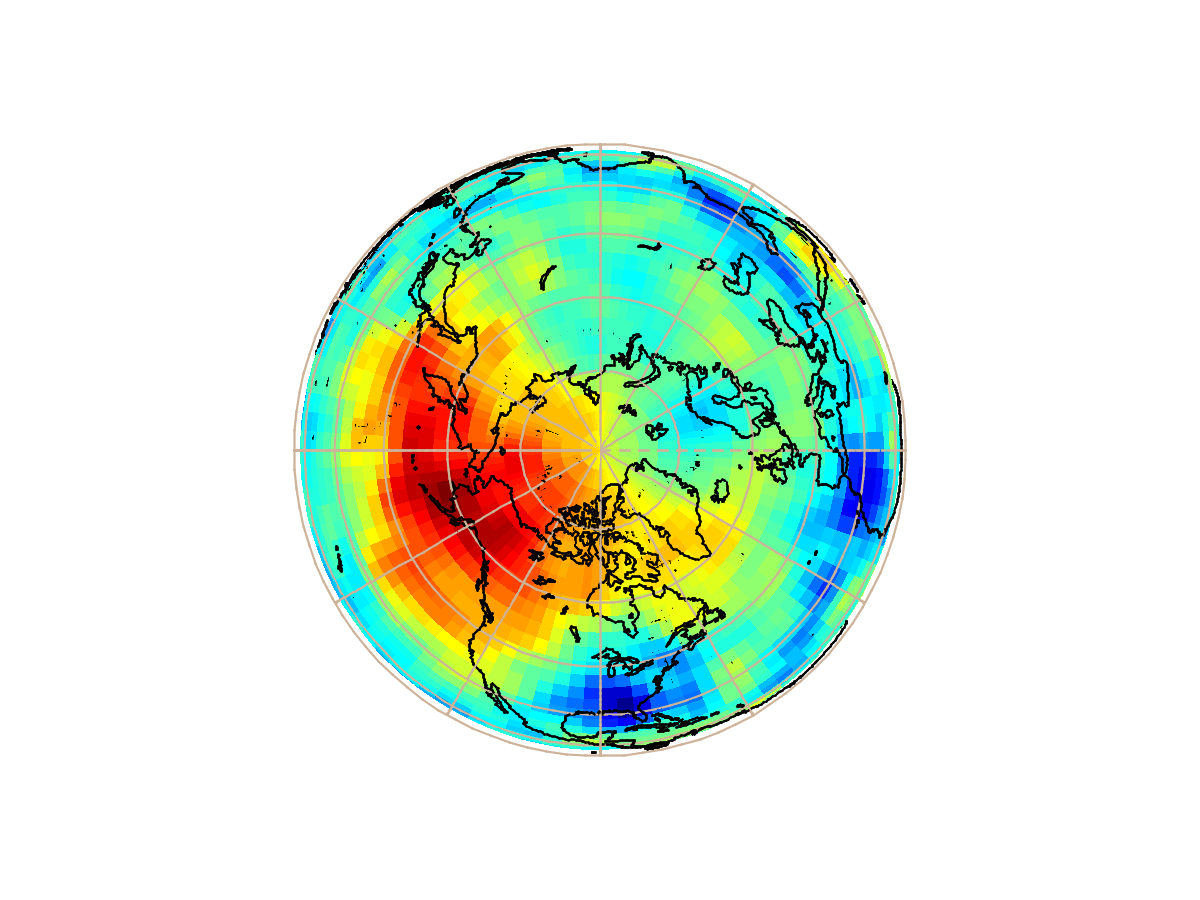}
  \caption{POD-4D-EnKF}
  \label{fig:T2}
\end{subfigure}
\begin{subfigure}{.5\textwidth}
  \centering
  \includegraphics[width=0.7\linewidth,height=0.7\linewidth]{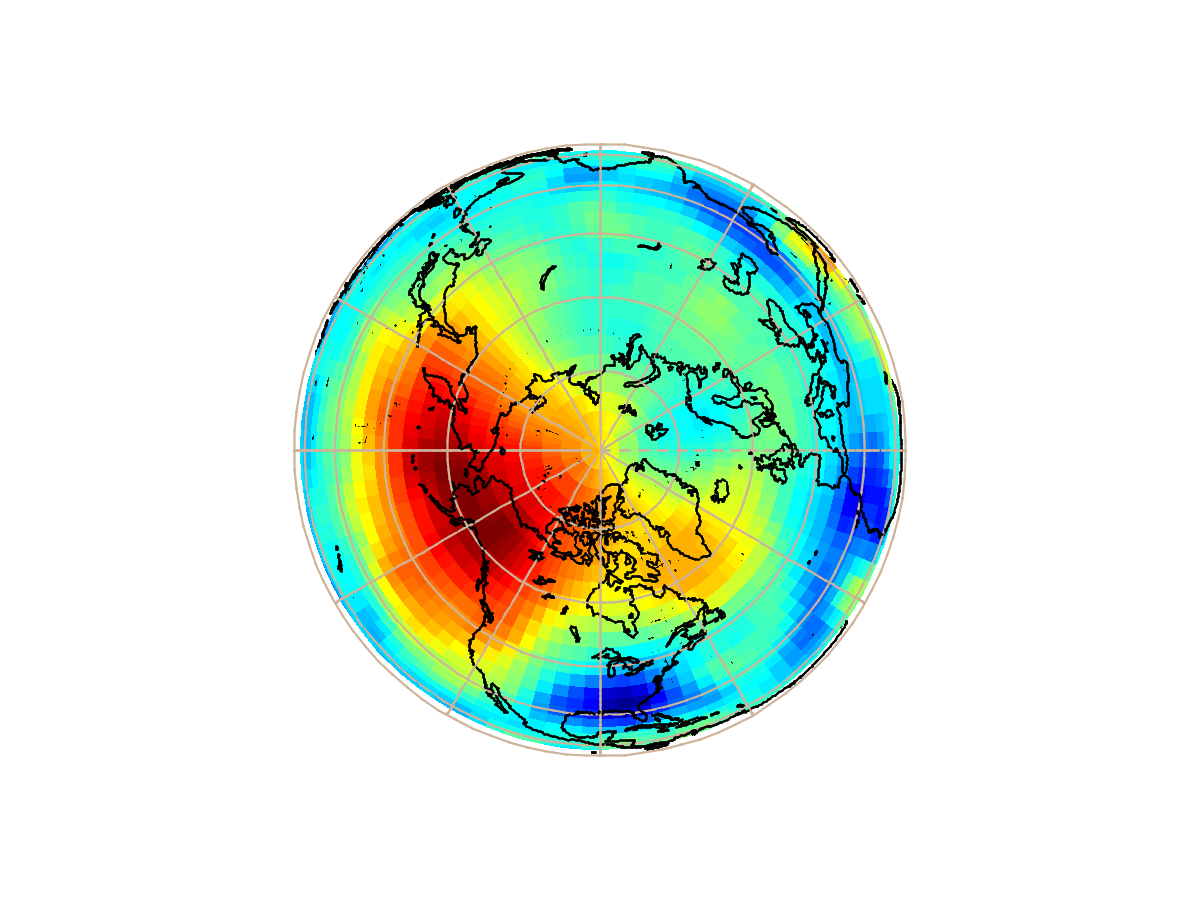}
   \caption{ISM}
  \label{fig:T3}
\end{subfigure}\hspace*{-7.8em}%
\begin{subfigure}{.5\textwidth}
  \centering
  \includegraphics[width=0.7\linewidth,height=0.7\linewidth]{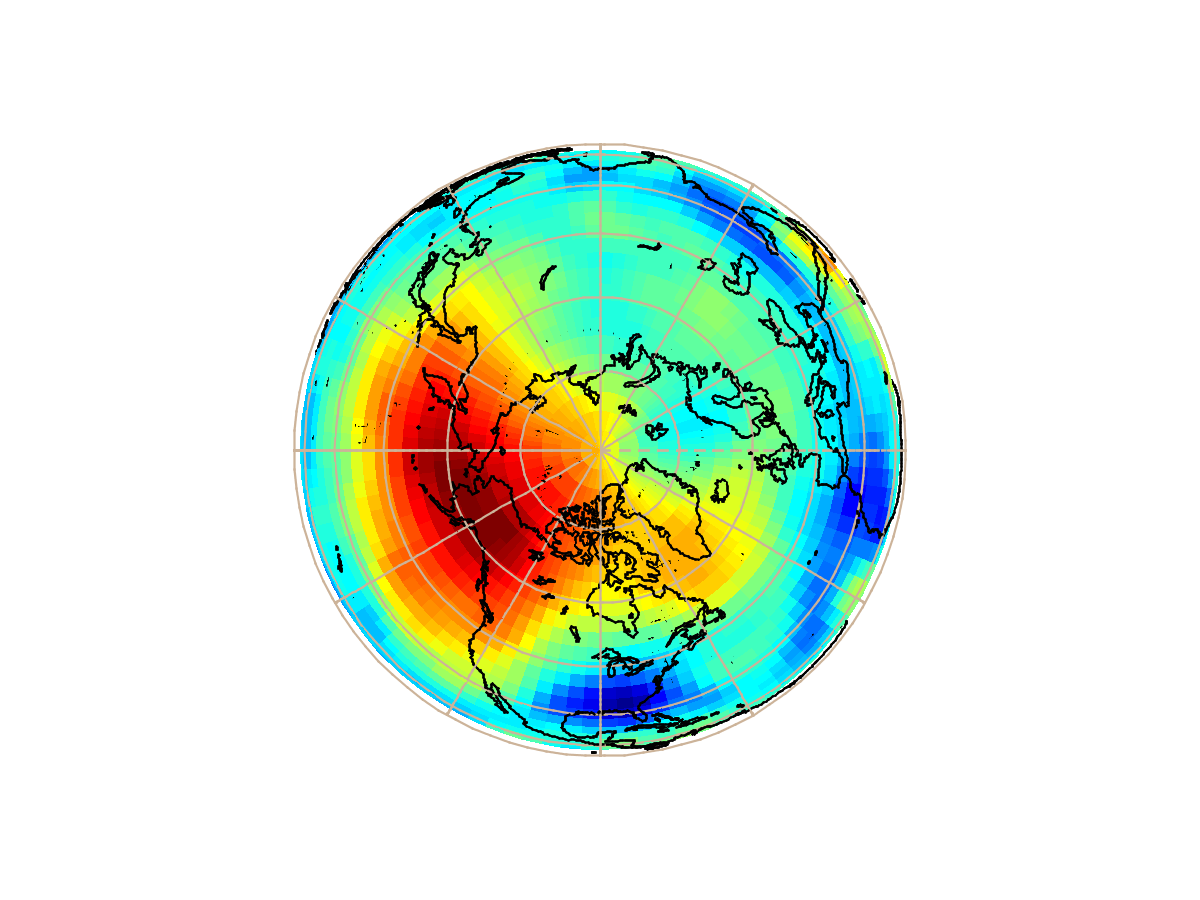}
   \caption{TR-4D-EnKF}
  \label{fig:T4}
\end{subfigure}
\caption{Estimated initial states ($\x^{a}_0$) for the temperature at the Earth's surface for the Background, POD-4D-EnKF, ISM and TR-4D-EnKF. After five iterations, the iterative methods ISM and 4D-TR-EnKF provide the most accurate results among the compared implementations.}
\label{fig:snapshots-initial-conditions-temperature}
\end{figure}

\begin{figure}[H]
\centering
\begin{subfigure}{.5\textwidth}
  \centering
  \includegraphics[width=0.8\linewidth,height=0.8\linewidth]{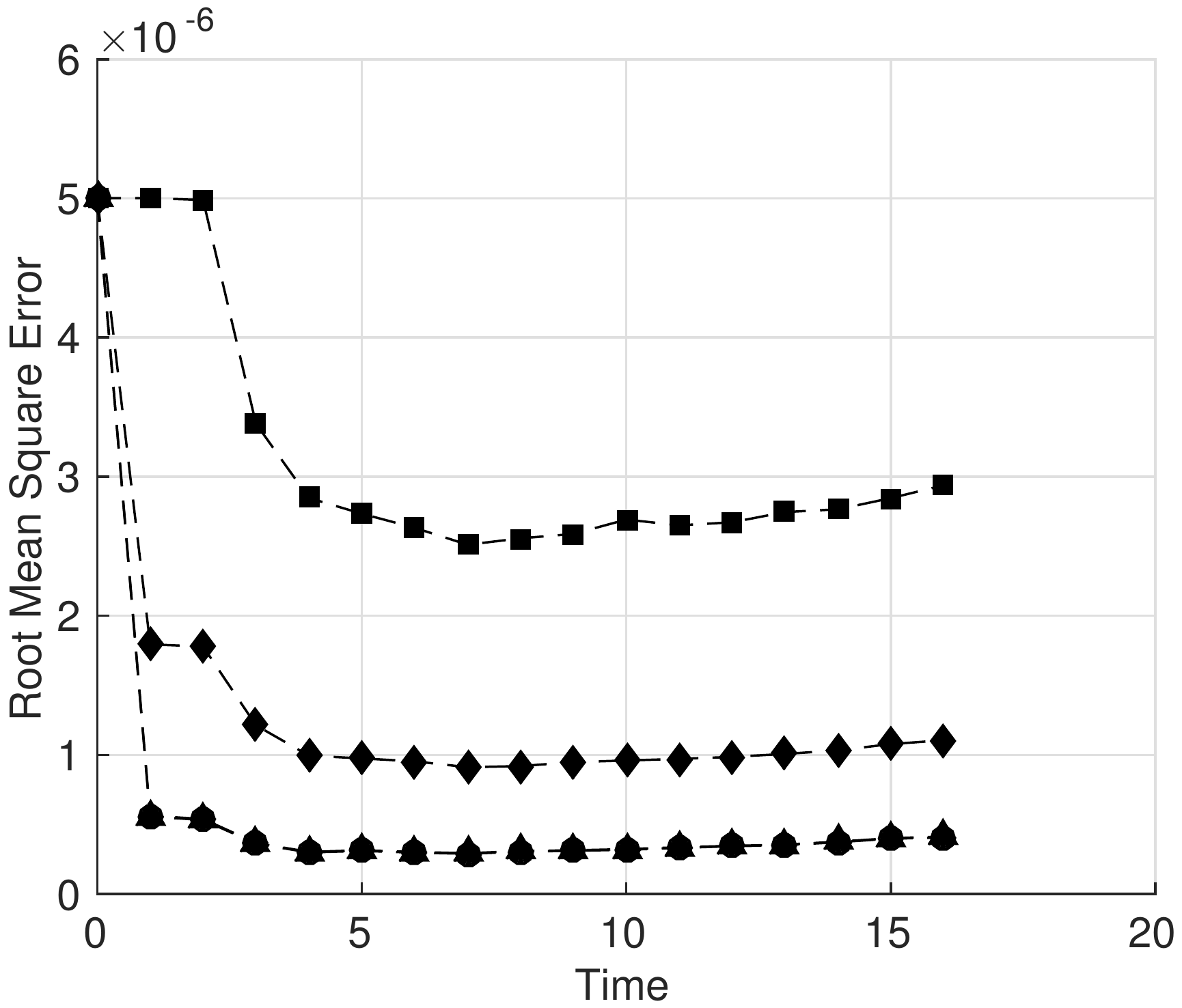}
  \caption{$\Nens = 10$}
  \label{fig:N10}
\end{subfigure}%
\begin{subfigure}{.5\textwidth}
  \centering
  \includegraphics[width=0.8\linewidth,height=0.8\linewidth]{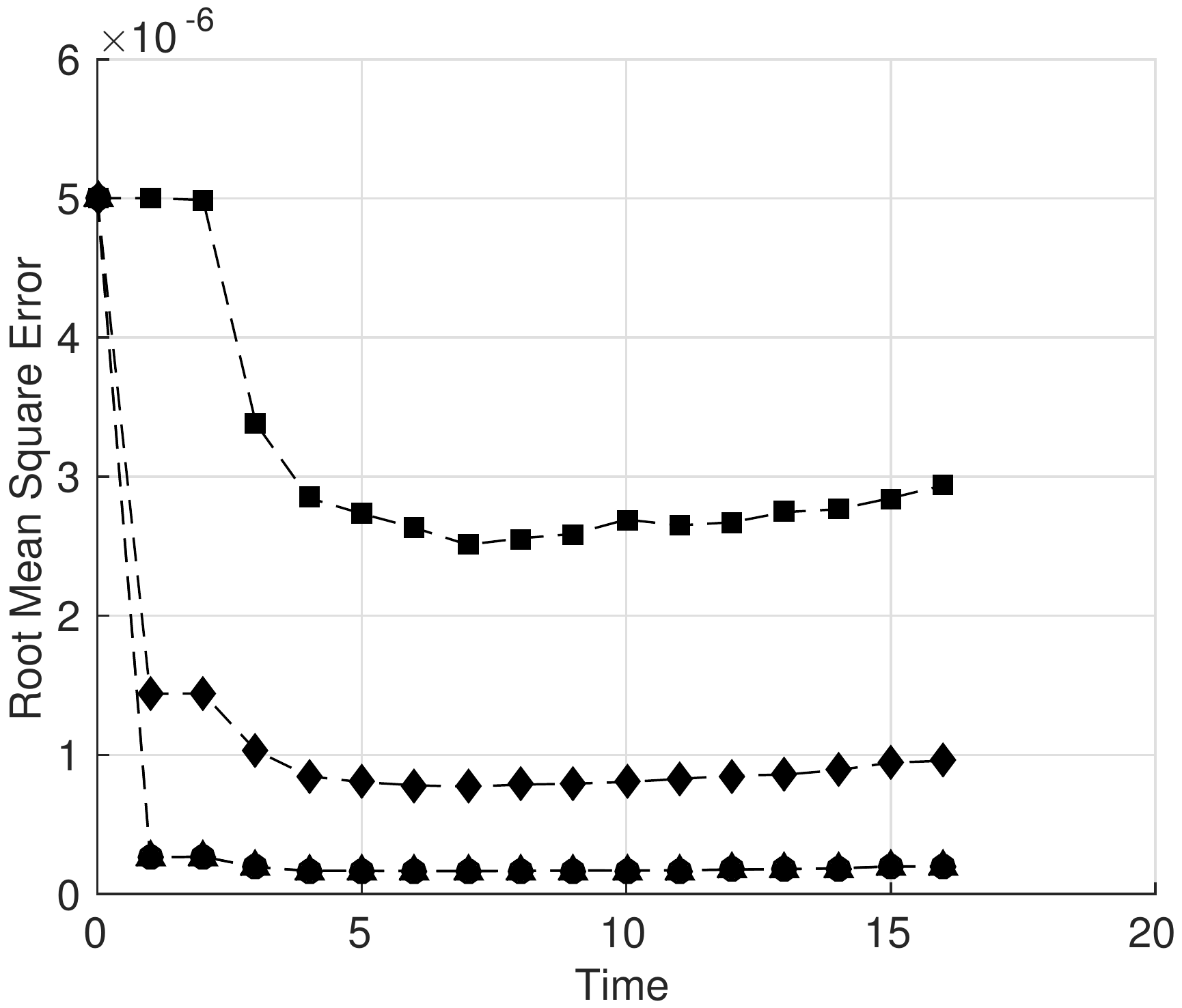}
  \caption{$\Nens=20$}
  \label{fig:N20}
\end{subfigure}
\begin{subfigure}{.5\textwidth}
  \centering
  \includegraphics[width=0.8\linewidth,height=0.8\linewidth]{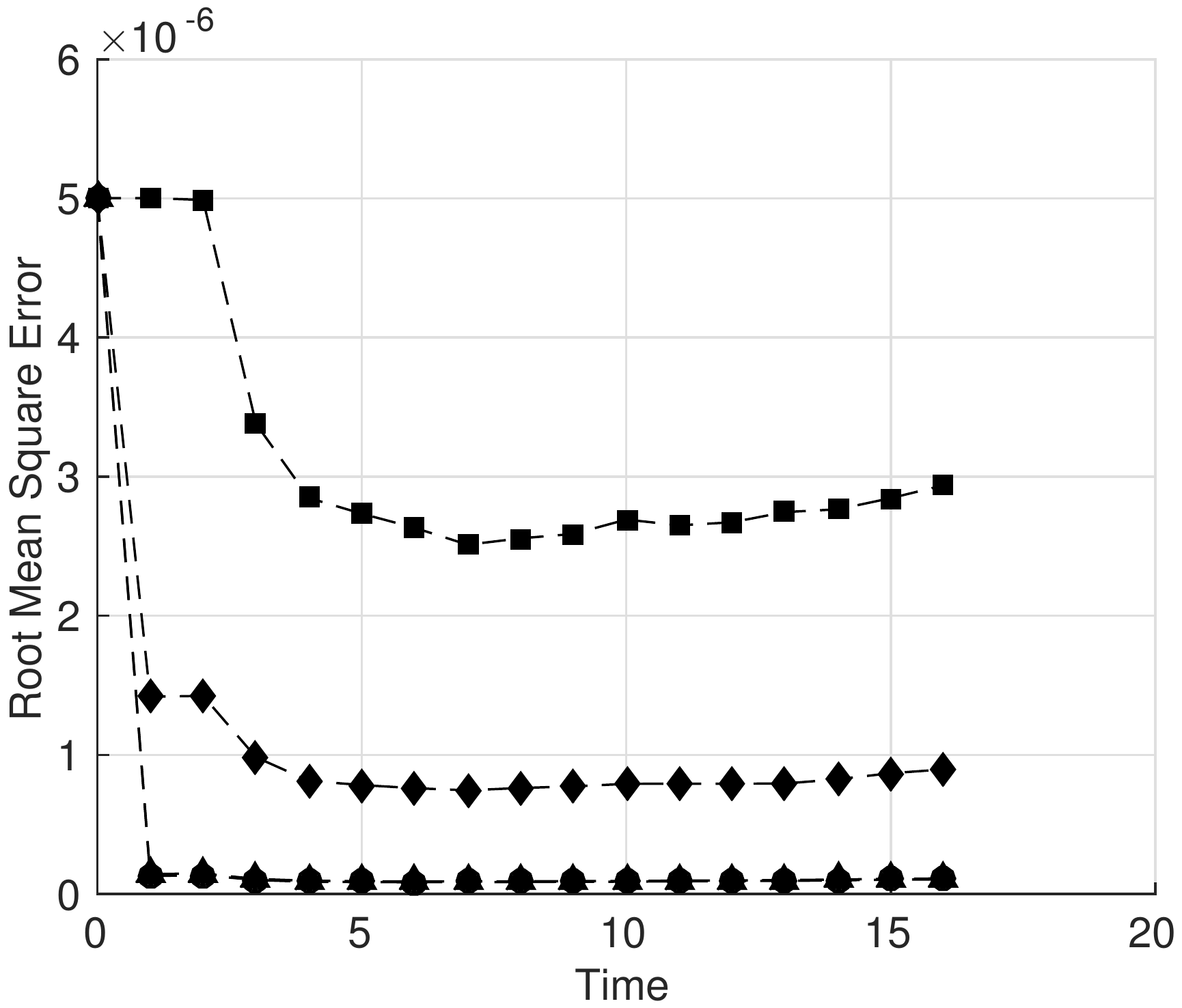}
   \caption{$\Nens=40$}
  \label{fig:N40}
\end{subfigure}%
\begin{subfigure}{.5\textwidth}
  \centering
  \includegraphics[width=0.8\linewidth,height=0.8\linewidth]{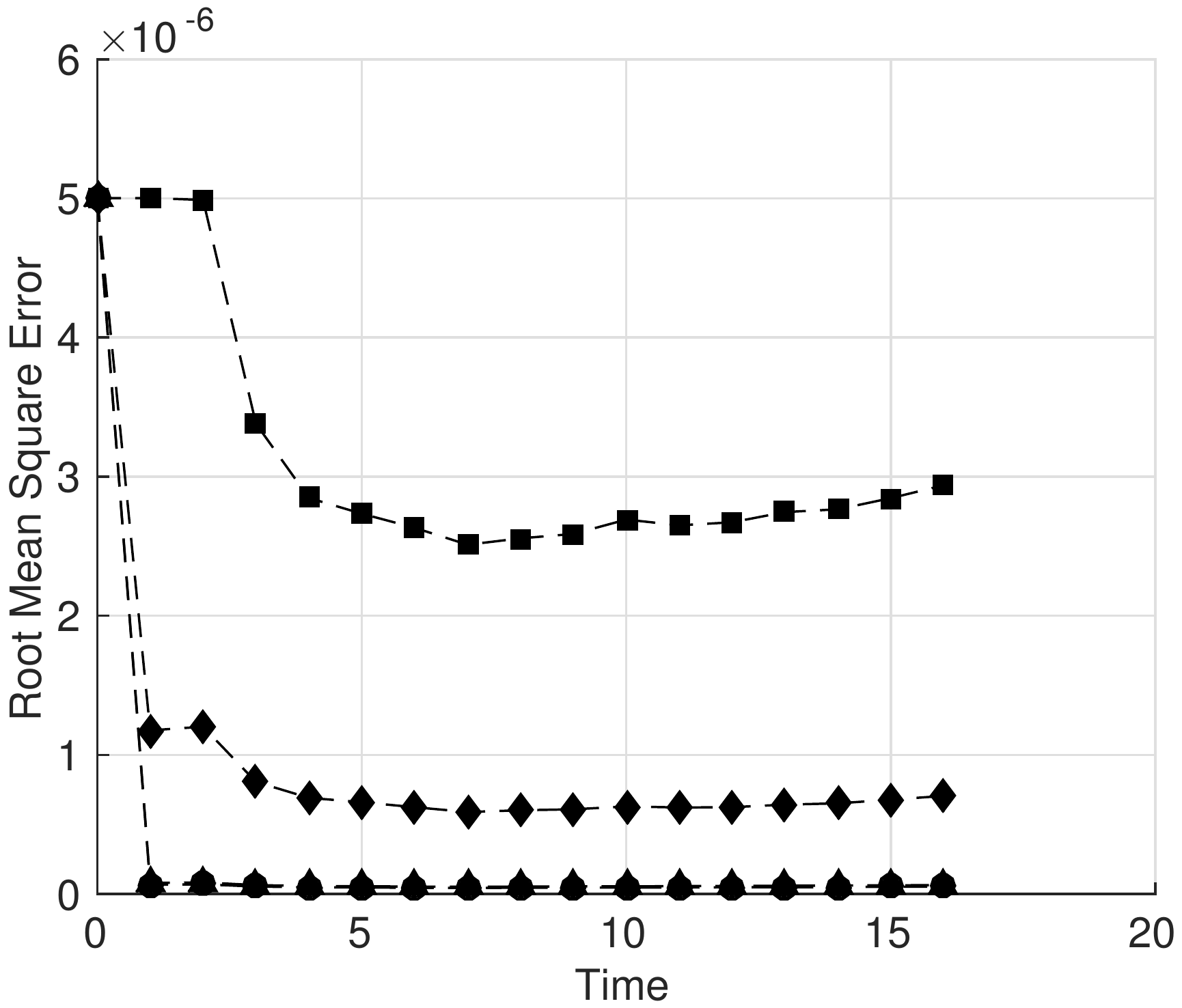}
   \caption{$\Nens=80$}
  \label{fig:N80}
\end{subfigure}
\caption{RMSE among observational times of the background state ($\blacksquare$), POD-4D-EnKF ($\blacklozenge$), ISM ($\blacktriangle$) and TR-4D-EnKF ($\bullet$) implementations for different ensemble sizes ($\Nens$). The most accurate results are obtained by the iterative methods. The analyses reported for the iterative methods are obtained after five iterations.}
\label{fig:vorticity-per-ensemble-size-time}
\end{figure}

\begin{figure}[H]
\centering
\begin{subfigure}{.5\textwidth}
  \centering
  \includegraphics[width=0.8\linewidth,height=0.8\linewidth]{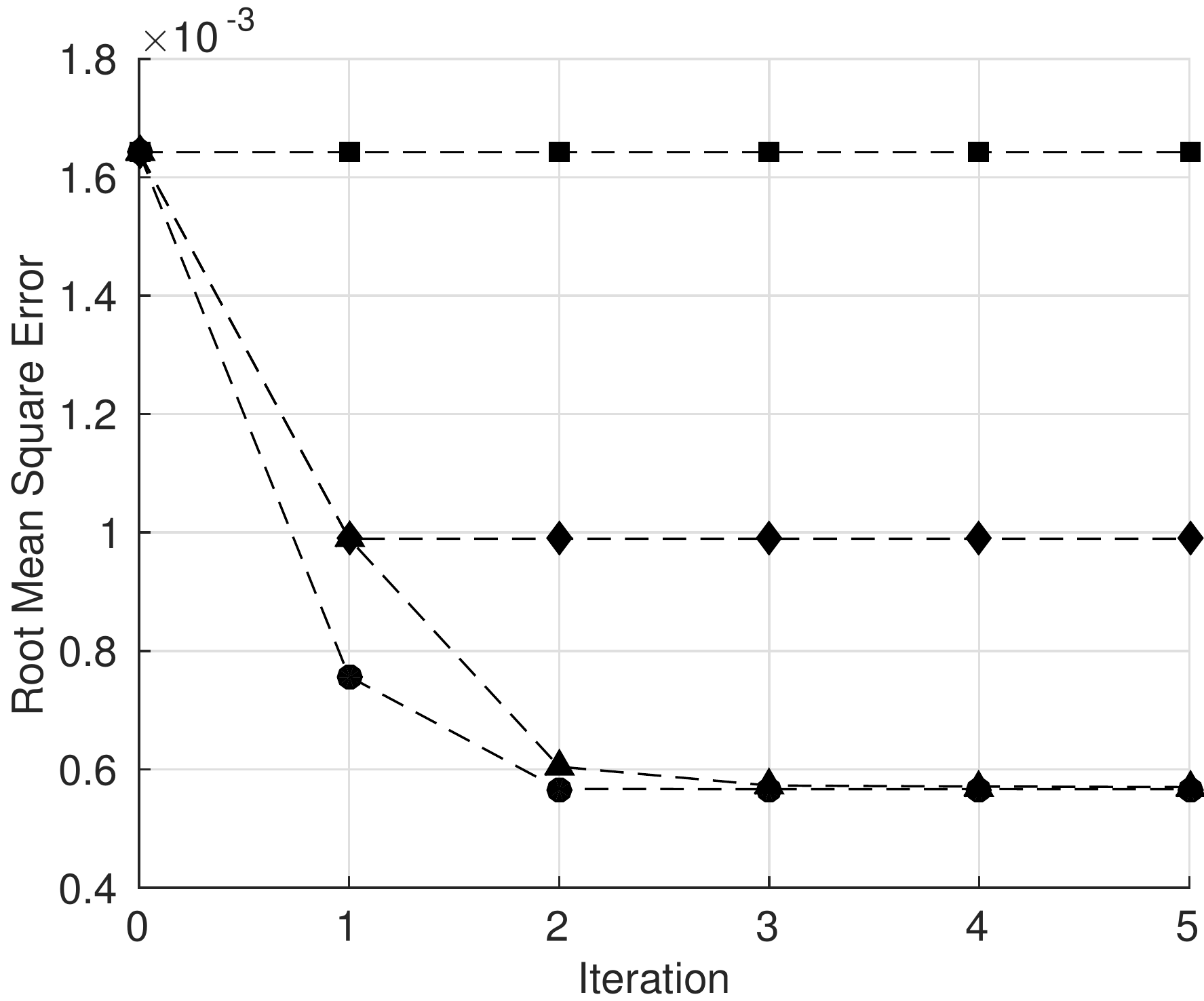}
  \caption{$\Nens = 10$}
  \label{fig:N10-O}
\end{subfigure}%
\begin{subfigure}{.5\textwidth}
  \centering
  \includegraphics[width=0.8\linewidth,height=0.8\linewidth]{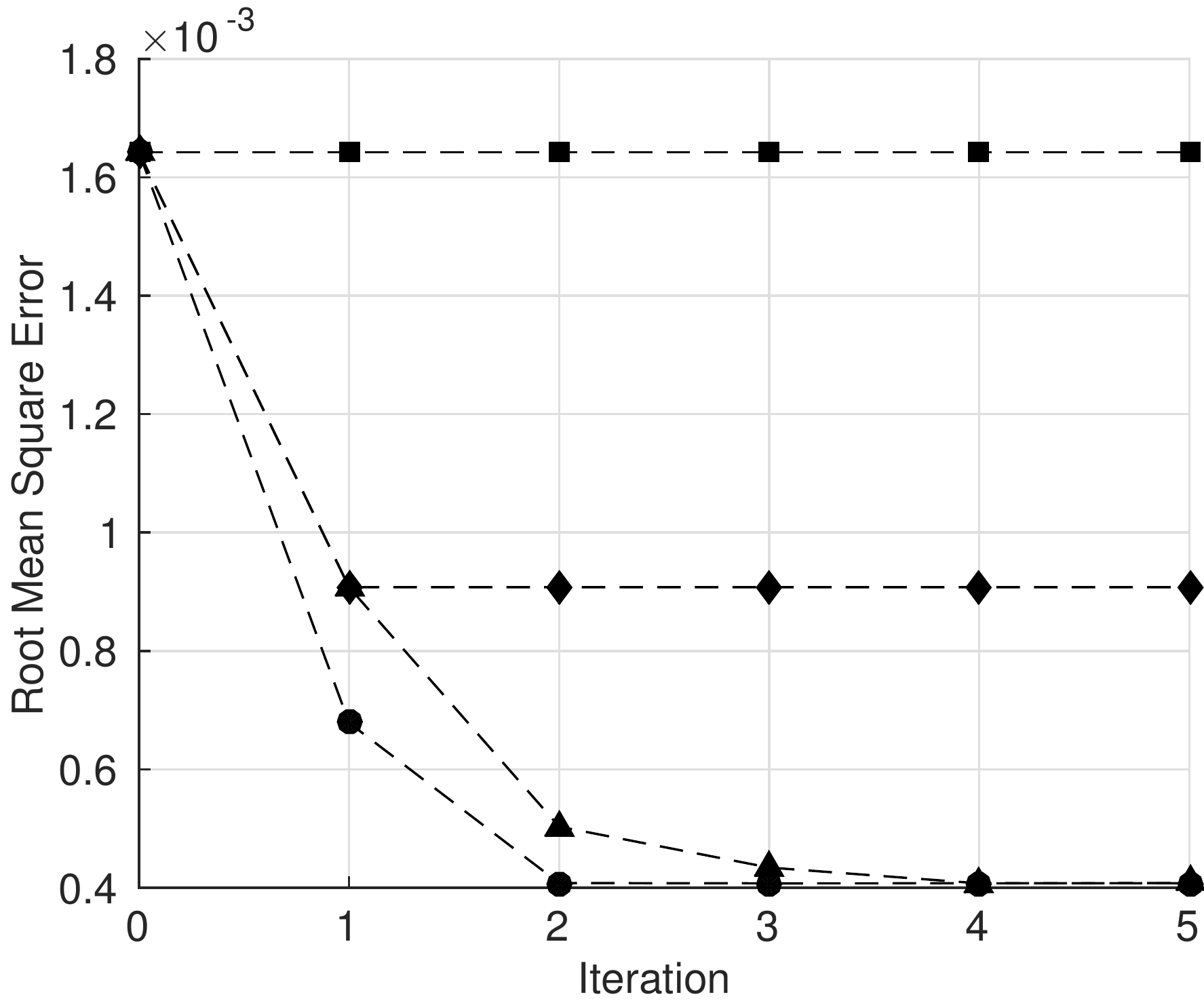}
  \caption{$\Nens=20$}
  \label{fig:N20-0}
\end{subfigure}
\begin{subfigure}{.5\textwidth}
  \centering
  \includegraphics[width=0.8\linewidth,height=0.8\linewidth]{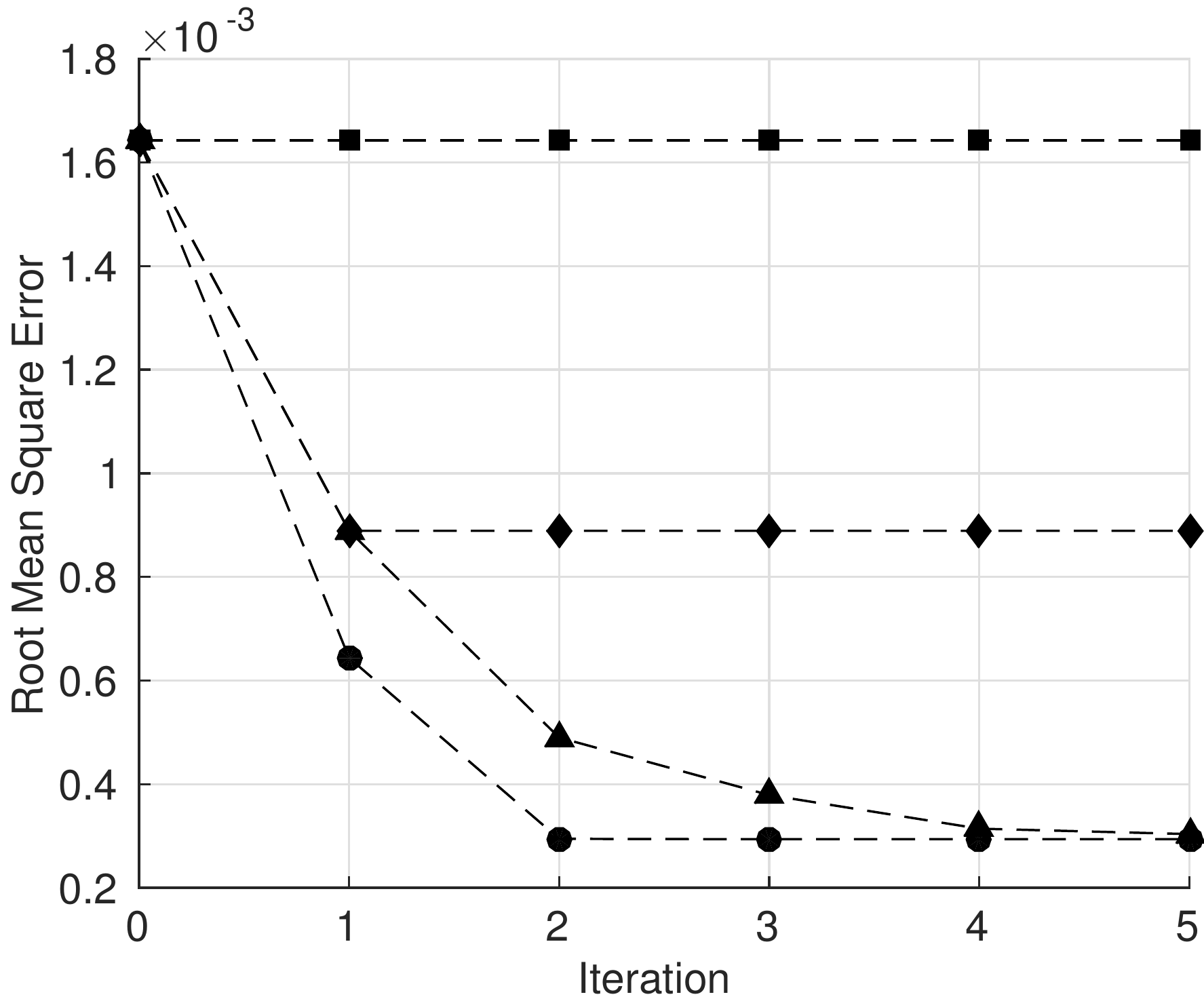}
   \caption{$\Nens=40$}
  \label{fig:N40-0}
\end{subfigure}%
\begin{subfigure}{.5\textwidth}
  \centering
  \includegraphics[width=0.8\linewidth,height=0.8\linewidth]{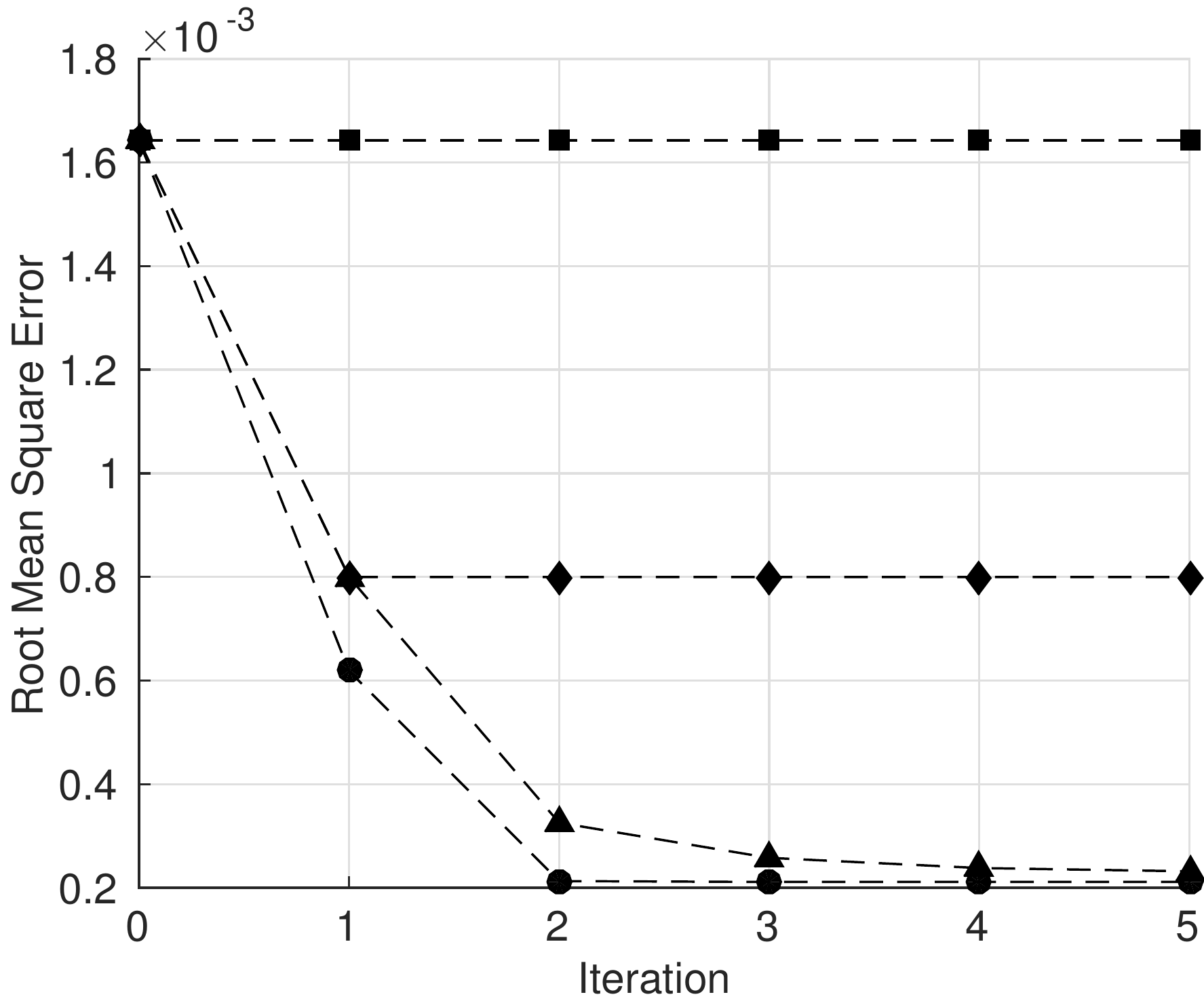}
   \caption{$\Nens=80$}
  \label{fig:N80-0}
\end{subfigure}
\caption{RMSE among iterations of the background state ($\blacksquare$), POD-4D-EnKF ($\blacklozenge$), ISM ($\blacktriangle$) and TR-4D-EnKF ($\bullet$) implementations for different ensemble sizes ($\Nens$). Since the POD-4D-EnKF is equivalent to one ISM iteration, its RMSE holds constant after the first iteration (for comparison purposes). The background is constant over all the iterations since it is the best estimation prior any measurement.}
\label{fig:vorticity-per-ensemble-size}
\end{figure}

\section{Conclusions}
\label{sec:conclusions}

This paper develops TR-4D-EnKF, an ensemble-based 4D-Var data assimilation method based on the trust region framework. The proposed implementation projects the model space onto the space spanned by the deviations of the ensemble members from the mean, as is typically done in 4D-EnKF implementations. A small optimization problem is solved in the ensemble space. At each iteration a new ensemble based surrogate model of the 4D-Var cost function is constructed, and the convergence is controlled by the trust region method. The trust region radius connects the optimal solution found in the ensemble space with the corresponding solution in the full model space. Moreover, the evolution of error statistics throughout iterations are captured by an empirical relation that uses the changes in trust region radius as a proxy for uncertainty decrease. Experimental results shows that the proposed implementation provide more accurate results than some of the best 4D-EnKF implementations available in the literature within a reasonable computational effort {and a lesser number of iterations}.

\section*{Acknowledgements}

This work was supported in part by awards NSF CCF--1218454, AFOSR FA9550--12--1--0293--DEF, AFOSR 12-2640-06, and by the Computational Science Laboratory at Virginia Tech. 
\appendix
\section{The trust region method}
\label{app:TR-method}
Consider the unconstrained minimization problem
\begin{eqnarray}
\label{eq:prel-function-tr}
\displaystyle 
\displaystyle \x^{*} = \underset{\x}{\arg\, \min}\,f(\x) \in \Re^{\Nstate \times 1}\,,
\end{eqnarray}
Trust Region based methods proceed as follows.
\begin{enumerate}
\item {\bf Initialization}. Define the initial solution $\x^{[0]} \in \Re^{\Nstate \times 1}$ and the parameters $\radius_{\rm max} \in (0,\, \infty)$ (maximum radius size), $\radius_{[0]} \in (0,\radius_{\rm max})$ (initial radius size), $\eta \in (0,2)$ (control variable for updating the solution), $0<\theta_1<\theta_2 < 1$ (control variables for updating the TR size), $\gamma_{\rm inc}$ (increasing factor of the TR size), $\gamma_{\rm dec}$ (decreasing factor of the TR size) and $j \leftarrow 0$ (iterate number).
\item {\bf Model generation}. Build the quadratic model $m_{[j]}(\dxt)$ as follows:
\begin{eqnarray}
\label{eq:prel-quadratic-model}
f \lp \x^{[j]}+\dxt \rp \approx m_{[j]}(\dxt) = f_{[j]}+\g_{[j]}^T \cdot \dxt + \frac{1}{2} \cdot \dxt^T \cdot \G_{[j]} \cdot \dxt \,.
\end{eqnarray}
In practice $f_{[j]} \approx f(\x^{[j]})$, $\g_{[j]} \approx \nabla f(\x^{[j]}) \in \Re^{m \times 1}$, and $\G_{[j]} \approx \nabla^2 f(\x^{[j]}) \in \Re^{m \times m}$, since the exact derivatives of \eqref{eq:prel-function-tr} are unavailable or difficult to compute.
\item {\bf Subproblem optimization}. Compute the optimal step via the solution of
\begin{subequations}
\begin{eqnarray}
\label{eq:prel-trust-region}
\dxt^{*} &=& \underset{\dxt} {\mathrm{arg \,min}} \,m_{[j]}(\dxt), \\
\label{eq:original-trust-region-constraint}
&& \textnormal{subject to} \quad \displaystyle \| \dxt \|\le \radius_{[j]}. 
\end{eqnarray}
\end{subequations}
\item {\bf Ratio computation}. Compute the ratio 
\begin{eqnarray}
\label{eq:prel-ratio}
\displaystyle 
\rho_{[j]} = \frac{f(\x^{[j]})-f(\x^{[j]}+\dxt^{*})}{m_{[j]}({\bf 0}_{\Nstate})-m_{[j]}(\dxt^{*})} \,,
\end{eqnarray}
where the numerator and denominator are often called the {\it actual} and {\it predicted} reduction.
\item {\bf Solution update}. Update the solution $\x^{(k)}$ according to
\begin{eqnarray}
\label{eq:solution-update}
\x^{[j+1]} &=& \begin{cases}
\x^{[j]} & \text{for $\rho \le \eta$} \\
\x^{[j]}+\dx^{*} & \text{otherwise} 
\end{cases} \,.
\end{eqnarray}
\item {\bf Radius update}. Update the radius $\radius_{[j]}$ according to 
\begin{eqnarray}
\label{eq:radius-update}
\radius_{[j+1]} &=& \begin{cases}
\radius_{[j]} \cdot \gamma_{\rm dec}  & \text{for $\rho < \theta_1$} \\
\radius_{[j]} & \text{for $\theta_1 \le \rho < \theta_2$} \\
\radius_{[j]} \cdot \gamma_{\rm inc} & \text{for $\theta_2 \le \rho \le 1$} \\
\end{cases} \,.
\end{eqnarray}
\item {\bf Iteration update}. Let $j \leftarrow j+1$ and go to 2. 
\end{enumerate}


\bibliographystyle{plain}

\end{document}